\documentclass{fic-l}
\usepackage{amsbsy,
   amssymb,amsmath,amscd}

\usepackage{latexsym}
\usepackage{fontenc}
\usepackage[mathscr]{euscript}
\usepackage[dvips]{epsfig}

\catcode`\@=11
\renewcommand{\over}{\@@over}
\renewcommand{\atop}{\@@atop}
\renewcommand{\above}{\@@above}
\renewcommand{\overwithdelims}{\@@overwithdelims}
\renewcommand{\atopwithdelims}{\@@atopwithdelims}
\renewcommand{\abovewithdelims}{\@@abovewithdelims}
\catcode`\@=13

 \font\thinlinefont=cmr5 \font\fiverm=cmr5
\input{pictex.tex}

\newcommand{\C}{\mathbb C}

\newcommand{\F}{\mathbb F}

\newcommand{\pP}{\mathbb P}
\newcommand{\Q}{\mathbb Q}
\newcommand{\R}{\mathbb R}
\newcommand{\Z}{\mathbb Z}

\newcommand{\ScrE}{{\mathscr{E}}}
\newcommand{\ScrF}{{\mathscr{F}}}
\newcommand{\ScrG}{{\mathscr{G}}}
\newcommand{\ScrO}{{\mathscr{O}}}

\newcommand{\ScrZ}{{\mathscr{Z}}}

\renewcommand{\mod}{\text{\rm mod}}

\theoremstyle{plain} {
   }

\newtheoremstyle{Theo}
{3pt} {3pt} {\itfamily} {} {\bf} {.} {.5em} {}

{}

\theoremstyle{remark} {
   }
{
   }
{
   \newtheorem*{Defn}{\bf Definition}}
\newtheorem*{Rems}{\bf Remarks}
\newtheorem*{Remark}{\bf Remark}

\newenvironment{Proof}{\bigskip\noindent\textrm{Proof:}}

\newcommand{\Le}{{{\mathchoice{%
   \,{\scriptstyle\le}\,}
   {\,{\scriptstyle\le}\,}%
   {\,{\scriptscriptstyle\le}\,}%
   {\,{\scriptscriptstyle\le}\,}}}}

\newcommand{\Ge}%
  {{\mathchoice{\,{\scriptstyle\ge}\,}
  {\,{\scriptstyle\ge}\,}%
  {\,{\scriptscriptstyle\ge}\,}%
  {\,{\scriptscriptstyle\ge}\,}}}
\newcommand{\w}%
  {{\mathchoice{\,{\scriptstyle\wedge}\,}
  {{\scriptstyle\wedge}}
  {{\scriptscriptstyle\wedge}}%
  {{\scriptscriptstyle\wedge}}}}
\newcommand{\mapright}[1]{\smash{%
  \mathop{\,\longrightarrow\,}\limits^{#1}}}

\newcommand{\lrdoublemaps}[4]%
{\raise3pt\hbox{$\mathop{\hbox to
#1pt{\rightarrowfill}\kern-35pt%
  \lower3.95pt\hbox to #2pt{\leftarrowfill}}%
  \limits_{#3}^{#4}$}}
\newcommand{\twoheaddownarrow}{%
  \raise8pt\hbox{$\vert$}
\raise1pt\hbox{\null\kern-4.65pt$\downarrow$}%
  \lower2pt\hbox{\null\kern-6pt$\downarrow$}}

\newcommand{\dashdown}{{\raise10pt\hbox{$%
   \vtop{\ialign{##\crcr
\hfil$\scriptscriptstyle\shortmid$\hfil\crcr
\noalign{\nointerlineskip}
    \hfil$\scriptscriptstyle\shortmid$\hfil\crcr
\noalign{\nointerlineskip}
    \hfil$\scriptscriptstyle\shortmid $\hfil\crcr
\noalign{\nointerlineskip}
  \hfil$\kern4.25pt\scriptscriptstyle\downarrow$
  \hfil\crcr}}$}\,}}

\newcommand{\dashup}{{\raise10pt\hbox{$%
   \vtop{\ialign{##\crcr
    \hfil$\scriptscriptstyle\uparrow$\hfil\crcr
\noalign{\nointerlineskip}
    \hfil$\scriptscriptstyle\shortmid$\hfil\crcr
\noalign{\nointerlineskip}
    \hfil$\scriptscriptstyle\shortmid $\hfil\crcr
\noalign{\nointerlineskip}
  \hfil$\kern4pt\scriptscriptstyle\shortmid$
  \hfil\crcr}}$}\,}}

\newcommand{\longmapright}[2]{\smash{\mathop
   {\hbox to #1pt{\rightarrowfill}}\limits^{#2}}}
\newcommand{\longmaprightsub}[2]%
  {\mathop{\hbox to #1pt{\rightarrowfill}}
     \limits_{#2}}

\newcommand{\Amalg}{\shortparallel\kern-7.80
true pt\lower.5 true pt%
  \hbox to 4.0 true pt{\hrulefill}\,\,}

\usepackage{graphicx}
\newcommand{\turndown}[1]{%
  \rotatebox[origin=c]{270}{\ensuremath#1}}
\newcommand{\hookdownarrow}{\turndown{\hookrightarrow}}
\newcommand{\turnup}[1]{%
  \rotatebox[origin=c]{90}{\ensuremath#1}}
\newcommand{\hookuparrow}{\turnup{\hookrightarrow}}

\def\lmapdown#1{\Big\downarrow%
  \llap{$\vcenter{%
   \hbox{$\scriptstyle#1\enspace$}}$}}
\def\rmapdown#1{\Big\downarrow%
   \rlap{$\kern-1.0pt\vcenter{%
   \hbox{$\scriptstyle#1$}}$}}

\newcommand{\doubleheaddownarrow}%
  {\hbox{$\Big\downarrow%
 \null\kern-7.00pt\lower2pt\hbox{$\downarrow$}$}}
\newcommand{\longtwoheadedrightarrow}[1]%
{\raise2.75pt\hbox to #1pt{\hrulefill}%
    \!\!\twoheadrightarrow}

\newcommand{\leftcapdownarrow}%
  {\raise9pt\hbox{$\scriptstyle\cap$}%
  \kern-4.00pt \Big\downarrow}
\newcommand{\rightcapdownarrow}%
  {\raise13pt\hbox{$\scriptscriptstyle\cap$}%
  \kern-8.30pt\Big\downarrow}

\newcommand{\Limright}[1]{{%
   \vtop{\ialign{##\crcr
    \hfil\rm lim\hfil\crcr
\noalign{\nointerlineskip\vskip2pt}
    \hfil\rightarrowfill\hfil\crcr
\noalign{\nointerlineskip\vskip1pt}
    \hfil$\scriptstyle #1$\hfil\crcr
\noalign{\nointerlineskip}\crcr}}}\,\,}

\newcommand{\Limleft}[1]{{%
   \vtop{\ialign{##\crcr
    \hfil\rm lim\hfil\crcr
\noalign{\nointerlineskip\vskip2pt}
   \hfil \leftarrowfill\hfil \crcr
\noalign{\nointerlineskip\vskip-1.50pt}
   \hfil$\scriptstyle#1$\hfil\crcr
   \noalign{\nointerlineskip}\crcr}}}}

\newcommand{\arrowsim}{%
   \smash{\mathop{\longrightarrow}%
  \limits^{\lower1.5pt\hbox{$\scriptstyle\sim$}}}}

\newcommand{\leftarrowsim}{%
\smash{\mathop{\longleftarrow}%
 \limits^{\lower1.5pt\hbox{$\scriptstyle\sim$}}}}

\newcommand{\nbar}{\bar{n}}
\newcommand{\ubar}{\bar{u}}
\newcommand{\vbar}{\bar{v}}

\newcommand{\Abar}{\bar{A}}
\newcommand{\Ebar}{\bar{E}}
\newcommand{\Fbar}{\bar{F}}
\newcommand{\Hbar}{\overline{H}}
\newcommand{\Pbar}{\bar{P}}
\newcommand{\Qbar}{\bar{Q}}

\newcommand{\Xbar}{\bar{X}}

\newcommand{\Res}{\text{\rm Res}}

\newcommand{\id}{\text{\rm id}}
\newcommand{\pr}{\text{\rm pr}}
\newcommand{\rat}{\text{\rm rat}}
\newcommand{\CH}{\text{\rm CH}}
\newcommand{\Alb}{\text{\rm Alb}}
\newcommand{\Ker}{\text{\rm Ker}}
\newcommand{\trans}{\text{\rm trans}}
\newcommand{\res}{\text{\rm res}}
\newcommand{\Char}{\text{\rm char}}
\newcommand{\cl}{\text{\rm cl}}
\newcommand{\et}{\text{\rm et}}
\newcommand{\Gal}{\text{\rm Gal}}

\newcommand{\cont}{\text{\rm cont}}
\newcommand{\Sym}{\text{\rm Sym}}
\newcommand{\Spec}{\text{\rm Spec}}
\newcommand{\Flat}{\text{\rm flat}}
\newcommand{\fl}{\text{\rm fl}}
\newcommand{\Pic}{\text{\rm Pic}}

\newcommand{\ST}{\text{\rm ST}}
\newcommand{\Br}{\text{\rm Br}}

\title{Local Galois Symbols on $E \times E$}
\author{Jacob Murre and Dinakar
Ramakrishnan
}
\date{}
\begin{document}
\begin{flushright}
{\it To Spencer Bloch, with admiration}
\end{flushright}
\begin{abstract}
This article studies the
Albanese kernel $T_F(E \times E)$, for an elliptic curve $E$ over a $p$-adic field $F$. The main result furnishes
information, for any odd prime $p$, about the kernel and image of the
Galois symbol map from $T_F(E \times E)/p$ to the Galois cohomology
group $H^2(F, E[p] \otimes E[p])$, for $E/F$ ordinary, without requiring that the $p$-torsion points are $F$-rational, or even that the Galois module $E[p]$ is semisimple.
A key step is to show that the image is zero when the finite Galois module $E[p]$ is acted on non-trivially by the pro-$p$-inertia group $I_p$. Non-trivial classes in the image are also constructed when $E[p]$ is suitably unramified.
A forthcoming sequel will deal with global questions.
\end{abstract}


\maketitle \overfullrule=5pt

\newcounter{mysubsection}[section]
\setcounter{mysubsection}{0}
\newcommand{\mysubsection}[1]{\refstepcounter{mysubsection} %
\vspace{0.5em} \noindent {\large\bf
\arabic{section}.\arabic{mysubsection} \hspace{0.4em} #1}
\par\vspace*{.2truecm}
}

\addcontentsline{toc}{section}{Introduction}
\section*{Introduction}

\bigskip

Let $E$ be an elliptic curve over a field $F$, $\overline F$ a
separable algebraic closure of $F$, and $\ell$ a prime different
from the characteristic of $F$. Denote by $E[\ell]$ the group of
$\ell$-division points of $E$ in $E(\overline F)$. To any
$F$-rational point $P$ in $E(F)$ one associates by Kummer theory a
class $[P]_\ell$ in the Galois cohomology group $H^1(F, E[\ell])$,
represented by the $1$-cocycle
$$
\beta_\ell: {\rm Gal}(\overline F/F) \rightarrow E[\ell], \, \,
\sigma \mapsto \sigma\left(\frac{P}{\ell}\right) - \frac{P}{\ell}.
$$
Here $\frac{P}{\ell}$ denotes any point in $E(\overline F)$ with
$\ell\left(\frac{P}{\ell}\right) = P$. Given a pair $(P,Q)$ of
$F$-rational points, one then has the cup product class
$$
[P,Q]_\ell : = \, [P]_\ell \cup [Q]_\ell \, \in \, H^2(F,
E[\ell]^{\otimes 2}).
$$

Any such pair $(P,Q)$ also defines a $F$-rational algebraic cycle
on the surface $E \times E$ given by
$$
\langle P,Q \rangle : = \, [(P,Q) - (P,0) - (0,Q) + (0,0)],
$$
where $[ \cdots ]$ denotes the class taken modulo {\it rational
equivalence}. It is clear that this {\it zero cycle of degree
zero} defines, by the parallelogram law, the trivial class in the
Albanese variety $Alb(E \times E)$. So $\langle P,Q\rangle$ lies
in the {\it Albanese kernel} $T_F(E \times E)$. It is a known, but
non-obvious, fact that the association $(P,Q) \to [P,Q]_\ell$
depends only on $\langle P,Q\rangle$, and thus results in the {\it
Galois symbol map}
$$
s_\ell: \, T_F(E \times E)/\ell \, \rightarrow \, H^2(F,
E[\ell]^{\otimes 2}).
$$
For the precise technical definition of $s_\ell$ we refer to Definition 2.5.3. In essence, $s_\ell$ is the restriction to the Albanese kernel $T_F(E\times E)$ of the cycle map over $F$ with target the continuous cohomology (in the sense of Jannsen \cite{J}) with $\Z/\ell$-coefficients. Due to the availability of the norm map for finite extensions (see Section 1.7), it suffices to study this map $s_\ell$ on the subgroup $ST_F(E\times E)$ of $T_F(E\times E)$ generated by symbols, whence the christening of $s_\ell$ as the {\it Galois symbol map}. Moreover, it is possible to study this map via a Kummer sequence (see Sections 2.2 and 2.3).

It is a conjecture of Somekawa and Kato that this map is always
injective (\cite{So}). It is easy to verify this when $F$ is $\C$ or
$\R$. In the latter case, the image of $s_\ell$ is non-trivial iff
$\ell = 2$ and all the $2$-torsion points are $\R$-rational, which
can be used to exhibit a non-trivial global $2$-torsion class in
$T_\Q(E \times E)$, whenever $E$ is defined by $y^2 =
(x-a)(x-b)(x-c)$ with $a,b,c \in \Q$. It should also be noted that
injectivity of the analog of $s_\ell$ fails for certain surfaces
occurring as quadric fibrations (cf. \cite{ParS}). However, the
general expectation is that such pathologies do not occur for {\it
abelian} surfaces.

\medskip

Let $F$ denote a non-archimedean local field, with ring of integers
$\ScrO_F$ and residual characteristic $p$. Let $E$ be a semistable
elliptic curve over $F$, and $\mathcal E[\ell]$ the kernel of
multiplication by $\ell$ on $\mathcal E$, which defines a finite
flat groupscheme $\ScrE[\ell]$ over $S={\rm Spec}(\ScrO_F)$. Let $F(E[\ell])$ denotes
the smallest Galois extension of $F$ over which all the
$\ell$-division points of $E$ are rational. It is easy to see that
the image of $s_\ell$ is zero if (i) $E$ has good reduction {\it
and} (ii) $\ell\ne p$, the reason being that the absolute Galois
group $G_F$ acts via its maximal unramified quotient
Gal$(F_{nr}/F)\simeq\hat{\Z}$, which has cohomological dimension
$1$. So we will concentrate on the more subtle $\ell=p$ case.

\medskip

\noindent{\bf Theorem A} \quad \it Let $F$ be a non-archimedean
local field of characteristic zero with residue field $\F_{q}$, $q=p^r$, $p$ odd.
Suppose $E/F$ is an elliptic curve over $F$, which has good,
ordinary reduction. Then the following hold:
\begin{enumerate}
\item[{(a)}] \, $s_p$ is injective, with image of $\F_p$-dimension $\leq
1$.
\item[{(b)}] \, The following are equivalent:
\begin{enumerate}
\item[{(bi)}] \, ${\text dim} \, {\text Im}(s_p)=1$
\item[{(bii)}] \, $F(E[p])$ is unramified over $F$, with the prime-to-$p$
part of $[F(E[p]):F]$ being $\leq 2$, and $\mu_p\subset F$.
\end{enumerate}
\item[{(c)}] \, If $E[p]\subset F$, then
$T_F(E \times E)/p \simeq \Z/p$ consists of symbols $\langle P, Q\rangle$
with $P, Q$ in $E(F)/p$.
\end{enumerate}
\rm

\medskip

Note that $[F(E[p]):F]$ is prime to $p$ iff the Gal$(\overline
F/F)$-representation $\rho_p$ on $E[p]$ is semisimple. We obtain:

\medskip

\noindent{\bf Corollary B} \, \it $T_F(E \times E)$ is $p$-divisible
when $E[p]$ is non-semisimple and wildly ramified. \rm

\medskip

When all the $p$-division points of $E$ are $F$-rational, the
injectivity part of part (a) has already been asserted, without
proof, in \cite{RS}, where the authors show the interesting result
that $T_F(E \times E)/p$ is a quotient of $K_2(F)/p$. Our techniques
are completely disjoint from theirs, and besides, the delicate part
of our proof of injectivity is exactly when $[F(E[p]):F]$ is
divisible by $p$, which is equivalent to the Galois module $E[p]$
being non-semisimple. Our results prove in fact that when $T_F(A)/p$
is non-zero, i.e., when $F(E[p])/F$ is unramified of degree $\leq
2$, it is isomorphic to the $p$-part Br$_F[p]$ of the Brauer group
of $F$, which is known (\cite{Tate3}) to be isomorphic to
$K_2(F)/p$.

\medskip

A completely analogous result concerning $s_\ell$ holds when $E/F$
has multiplicative reduction, for both $\ell=p$ and $\ell\ne p$, and
this can be proved by arguments similar to the ones we use in the
ordinary case. However, there is already an elegant paper of Yamazaki
(\cite{Y}) giving the analogous result (in the multiplicative
case) by a different method, and we content ourselves to a very brief discussion in
section 7 on how to deduce this analogue from \cite{Y}.

\medskip

The relevant preliminary material for the paper is assembled in the
first two sections and in the Appendix. We have, primarily (but not
totally) for the convenience of the reader, supplied proofs of
various statements for which we could not find published references,
even if they are apparently known to experts or in the {\it
folklore}. After this the proof of Theorem A is achieved in the following sections
via a number of Propositions. In section 3, we first prove, via the key
Proposition C, that the image of $s_p$ is at most one and establish part of the
implication $(bi) \implies (bii)$; the remaining part of this implication is
proved in section 4 via Proposition D. The converse implication $(bii)
\implies (bi)$ and part $(c)$ are proved in section 5. Finally, part $(a)$,
giving the injectivity of $s_p$, is
proved in section 6.

\medskip

In a sequel we will use these results in conjunction with others
(including a treatment of the case of supersingular reduction,
$p$-adic approximation, and a local-global lemma) and prove two
global theorems about the Galois symbols on $E\times E$ modulo $p$,
for any odd prime $p$.

\medskip

The second author would like to acknowledge support from the
National Science Foundation through the grants DMS-0402044 and
DMS-$0701089$, as well as the hospitality of the {\it DePrima
Mathematics House} at Sea Ranch, CA. The first author would like to
thank Caltech for the invitation to visit Pasadena, CA, for several
visits during the past eight years. This paper was finalized in May
2008. We are grateful to various mathematicians who pointed out errors
in a previous version, which we have fixed. Finally, it is a pleasure to thank
the referee of this paper, and also Yamazaki, for their thorough reading and for pointing out the need
for further explanations; one such comment led to our noticing a small gap,
which we have luckily been able to take care of with a modified argument.

\medskip

It gives us great pleasure to dedicate this article to Spencer
Bloch, whose work has inspired us over the years, as it has so many
others. The first author (J.M.) has been a friend and a fellow
{\it cyclist} of Bloch for many years, while the second author (D.R.) was
Spencer's post-doctoral mentee at the University of Chicago during
1980-82, a period which he remembers with pleasure.

\bigskip

\setcounter{section}{0}
\section{Preliminaries on Symbols}
\label{sec0}

\mysubsection{Cycles} \label{subsec1.1}

\medskip

Let $F$ be a field and $X$ a smooth projective variety of dimension
$d$ and defined over $F$. Let $0\Le i\Le d$ and $\ScrZ_i(X)$ be the
group of algebraic cycles on $X$ defined over $F$ and of dimension
$i$, i.e. the free group generated by the subvarieties of $X$ of
dimension $i$ and defined over $F$ (see \cite{F}, section 1.2). A cycle of
dimension $i$ is called rationally equivalent to zero over $F$ if
there exists $T\in\ScrZ_{i+1}(\pP^1\times X)$ and two points $P$
and $Q$ on $\pP^1$, each rational over $F$, such that $T(P)=Z$ and
$T(Q)=0$, where for $R\in\pP^1$ we put $T(R)=\pr_2(T.(R\times X))$
in the usual sense of calculation with cycles (see \cite{F}). Let
$\ScrZ_i^{\rat}(X)$ be the subgroup of $Z_i(X)$ generated by the
cycles rationally equivalent to zero and $\CH_i(X)=\ScrZ_i(X)/
\ScrZ_i^{\rat}(X)$ the Chow group of rational equivalence classes of
$i$-dimensional cycles on $X$ defined on $F$.

If $K\supset F$ is an extension then we write $X_K=X\times_F K$,
$Z_i(X_K)$ and $\CH_i(X_K)$ for the corresponding groups defined
over $K$ and if $\Fbar$ is the algebraic closure of $F$ then we
write $\Xbar=X_{\Fbar}$, $\ScrZ_i(\Xbar)$ and $\CH_i(\Xbar)$.

\medskip

\mysubsection{The Albanese kernel} \label{subsec1.2}

\medskip

In this paper we are only concerned with the case when $X$ is of
dimension $2$, i.e. $X$ is a {\em surface} (and in fact we shall
only consider very special ones, see below) and with groups of
$0$-dimensional cycles contained in $\CH_0(X)$. In that case note
that if $Z\in\ScrZ_0(X)$ then we can write $Z=\sum n_i P_i$ when
$P_i\in X(\Fbar)$, i.e. the $P_i$ are points on $X$ -- but
themselves only defined (in general) over an extension field $K$ of
$F$. Put $A_0(X)=\Ker(\CH_0(X){\overset\deg\longrightarrow} \Z)$
where $\deg(Z)=\sum n_i \deg(P_i)$, i.e. $A_0(X)$ are the cycle
classes of degree zero and put further $T(X):=\Ker(A_0(X)\to
\Alb(X))$, where $\Alb(X)$ is the Albanese variety of $X$ (see
\cite{Bl}). Again if $K\supset F$, write $A_0(X_K)$, $T(X_K)$ or
also $T_K(X)$ for the corresponding groups over $K$.

It is known that both $A_0(\Xbar)$ and $T(\Xbar)$ are {\em
divisible} groups (\cite{Bl}). Moreover if the group of
transcendental cycles $H^2(\Xbar)_{\trans}$ is non-zero, and if
$\Fbar$ is a universal domain (ex. $\Fbar=\C$), the group
$T(\Xbar)$ is ``huge'' (infinite dimensional) by Mumford (in
characteristic zero) and Bloch (in general) (\cite{Bl}, 1.22).

\bigskip

\mysubsection{Symbols} \label{subsec1.3}

\medskip

We shall only be concerned with surfaces $X$ which are abelian, even
of the form $A=E\times E$, where $E$ is an {\em elliptic curve}
defined over $F$. If $P,Q\in E(F)$, put
$$
\left<P,Q\right>:=(P,Q)-(P,O)-(O,Q)+(O,O)
$$
where $O$ is the origin on $E$, and the addition is the {\em
addition of cycles} (or better: cycle classes). Clearly
$\left<P,Q\right>\in T_F(A)$ and $\left<P,Q\right>$ is called the
{\em symbol} of $P$ and $Q$.

\begin{Defn}
$ST_F(A)$ denotes the subgroup of $T_F(A)$ generated by symbols
$\left<P,Q\right>$ with $P\in E(F)$, $Q\in E(F)$. It is called the
{\em symbol group} of $A=E\times E$.
\end{Defn}

\mysubsection{Properties and remarks} \label{subsec1.4}

\begin{description}
\item[1.]
The symbol is {\em bilinear} in $P$ and $Q$. For instance
$$
\left<P_1+P_2,Q\right>=\left<P_1,Q\right>+ \left<P_2,Q\right>
$$
(where now the first $+$ sign is addition on $E$!).

\item[2.]
The symbol is the Pontryagin product
$$
\left<P,Q\right>= \{(P,0)-(0,0)\}*\{(0,Q)-(0,0)\}
$$

\item[3.]
If $F$ is finitely generated over the prime field, then it follows
from the theorem of Mordell-Weil that $ST_F(A)$ is finitely
generated.

\item[4.]
Note the similarity with the group $K_2(F)$ of a field. Also the
symbol group $ST_F(A)$ is related to, but different from, the group
$K_2(E\times E)$ defined by Somekawa (\cite{So}); compare also with
\cite{RS}.
\end{description}

\medskip

\mysubsection{Restriction and norm/corestriction} \label{subsec1.5}

\medskip

Let $K\supset F$. Consider the morphism $\varphi_{K/F}\colon\,
A_K\to A_F$.

\noindent a.\enspace This induces homomorphisms, called the {\em
restriction} homomorphisms:
$$
\res_{K/F}=\varphi_{K/F}^*\colon\CH_i (A_F)\to \CH_i(A_K), T_F(A)\to
T_K(A)
$$
and $ST_F(A)\to ST_K(A)$ (see \cite{F}, section 1.7).

\noindent b.\enspace Also, when $[K:F] < \infty$, we have the {\em norm} or {\em
corestriction} homomorphisms
$$
N_{K/F}=(\varphi_{K/F})_*\colon\,CH_i(A_K)\to CH_i(A_F)\text{ and }
T(A_K)\to T_F(A)
$$
(see \cite{F}, section 1.4, page 11 and 12).

\begin{Rems}

\medskip

\begin{description}
\item[1.]
If $[K\colon\,F]=n$ we have $\varphi_*\circ\varphi^*=n$

\item[2.]
It is {\em not} clear if $N_{K/F}=\varphi_*$ induces a homomorphism
on the {\em symbol groups themselves}. Note however that there is
such norm map for cohomology (\cite{Serre3}, p. 127) and for
$K_2$-theory (\cite{M}, p. 137).
\end{description}
\end{Rems}

\medskip

\mysubsection{Some preliminary lemmas} \label{subsec1.6}

\medskip

\noindent{\bf Lemma 1.6.1} \, \it Let $Z\in T_F(A)$. Suppose we can
write $Z=\sum\limits_{i=1}^N
(P_i,Q_i)-\sum\limits_{i=1}^N(P'_i,Q'_i)$ with $P_i,P'_i,Q_i$ and
$Q'_i$ all in $E(F)$ for $i=1,\dotsc,N$. Then $Z\in ST_F(A)$. \rm

\medskip

\begin{Proof}
Since $Z\in T_F(A)$ we have that
$\sum\limits_{i=0}^{N}P_i=\sum\limits_{i=0}^{N} P'_i$ and
$\sum\limits_{i=1}^{N}Q_i=\sum\limits_{i=1}^{N} Q'_i$ as {\em sum of
points on $E$}. From this it follows immediately that we can rewrite
$Z=\sum\limits_{i=1}^{N}\{(P_i,Q_i)-(P_i,
0)-(0,Q_i)+(0,0)\}-\sum\limits_{i=1}^{N}
\{(P'_i,Q'_i)-(P'_i,0)-(0,Q'_i)+(0,0)
\}=\sum\limits_{i=1}^{N}\left<P_i,Q_i\right>-
\sum\limits_{i=1}^{N}\left<P'_i,Q'_i\right>$; hence $Z\in ST_F(A)$.
\end{Proof}

\noindent{\bf Corollary 1.6.2} \, \it Over the algebraic closure we
have
$$
T(\Abar) =\Limright{K} \{\res_{\Fbar/K}ST_K(A)\}
$$
where the limit is over all finite extensions $K\supset F$. \rm

\medskip

\begin{Proof}
Immediate from Lemma 1.6.1.
\end{Proof}

\medskip

\setcounter{subsection}{6} \mysubsection{Norms of symbols}
\label{subsec1.7}

\medskip

Next we turn our attention again to $T(A)=T_F(A)$ itself. If
$K\supset F$ is a finite extension and $P,Q\in E(K)$ then consider
$N_{K/F}\left(\left< P,Q\right>\right)\in T_F(A)$. Let
$ST_{K/F}(A)\subset T_F(A)$ be the subgroup of $T_F(A)$ generated by
such elements (i.e., coming as norms of the symbols from finite
field extensions $K\supset F$). Note that clearly
$ST_{F/F}(A)=ST_F(A)$ and also that $ST_{K/F}(A)$ consists of the
norms of elements of $ST_K(A)$.

\medskip\noindent
{\bf Lemma 1.7.1} {\it With the above notations let $T'_F(A)$ be the
subgroup of $T_F(A)$ generated by all the subgroups of type
$ST_{K/F}(A)$ of $T_F(A)$ for $K\supset F$ finite (with
$K\subset\Fbar$). Then $T'_F(A)=T_F(A)$.}

\medskip

\begin{Proof}
For simplicity we shall (first) assume $\Char(F)=0$. If $\cl(Z)\in
T_F(A)$, then $\cl(Z)$ is the (rational equivalence) class of a
cycle $Z\in\ScrZ_0(A_F)$ and moreover $Z=Z'-Z''$ with $Z'$ and $Z''$
positive (i.e. ``effective''), of the same degree and both
$Z'\in\ScrZ_0(A_F)$ and $Z''\in\ScrZ_0(A_F)$. Fixing our attention
on $Z'\in\ScrZ_0(A_F)$ we can, by definition, write
$Z'=\sum\limits_{\alpha}Z'_\alpha$ where the $Z'_\alpha$ are
($0$-dimensional) subvarieties of $A$ and {\em irreducible} over $F$
(Remark: in the terminology of Weil's Foundations \cite{W} the $Z'$
is a ``rational chain'' over $F$ and the $Z'_\alpha$ are the ``prime
rational'' parts of it, see \cite{W}, p. 207). For each $\alpha$ we
have $Z'_\alpha=\sum\limits_{i=1}^{n_\alpha}(P'_{\alpha
i},Q'_{\alpha i})$ where the $(P'_{\alpha i},Q'_{\alpha i})\in
A(\Fbar)$ is a set of points which form a {\em complete set of
conjugates over} $F$ (see \cite{W}, 207). Note that also
$\sum\limits_{i=1}^{n_\alpha}(P'_{\alpha i},0)\in\ScrZ_0(A_{F})$ and
similarly $\sum\limits_{i=1}^{n_\alpha}(0,Q'_{\alpha i})
\in\ScrZ_0(A_F)$ (Note: these cycles are rational, but not
necessarily prime rational.) For each $\alpha$ fix an arbitrary
$i(\alpha)$ and put $K_{\alpha i(\alpha)}= F(P'_{\alpha
i(\alpha)},Q'_{\alpha i(\alpha)})$ and consider now the {\em cycle}
$(P'_{\alpha i(\alpha)},Q'_{\alpha
i(\alpha)})\in\ScrZ_0(A_{K_{\alpha i(\alpha)}})$. We have by
definition (see \cite{F}, section 1.4, p. 11),
$Z'_\alpha=N_{K_{\alpha i(\alpha)}/F} (P'_{\alpha
i(\alpha)},Q'_{\alpha i(\alpha)})$. Furthermore if we put
$$
Z^\ast_{\alpha i(\alpha)}=\left<P'_{\alpha i(\alpha)},Q'_{\alpha
i(\alpha)}\right>\in ST_{K_{\alpha i(\alpha)}}(A)
$$
then in $T_F(A)$ we have
$$
N_{K_{\alpha i(\alpha)}/F}(Z^\ast_{\alpha
i(\alpha)})=\sum\limits_{i=1}^{n_\alpha} \left<P'_{\alpha
i},Q'_{\alpha i}\right>
$$
(again by the definition of the $N_{-/F}$) and clearly the cycle is
in $ST_{K_{\alpha i(\alpha)}/F}(A)$, i.e. in $T'_F(A)$. Now doing
this for every $\alpha$ and treating similarly
$Z''=\sum\limits_{\beta}Z''_\beta$, we have, since $Z\in T_F(A)$,
that
$$
Z=Z'-Z''=\sum\limits_{\alpha}Z'_\alpha- \sum\limits_{\beta}
Z'_\beta=\sum\limits_{\alpha}\sum\limits_{i} \left<P'_{\alpha
i},Q'_{\alpha i}\right>- \sum\limits_{\beta}\sum\limits_{j}\left<
P''_{\beta j},Q''_{\beta j}\right>.
$$
Hence $Z\in T'_F(A)$, which completes the proof (in $\Char\,0$).
\end{Proof}

\qed

\medskip

\begin{Remark}
If $\Char(F)=p>0$, then we have that
$Z'_\alpha=p^{m_\alpha}\sum\limits_{i} (P'_{\alpha i},Q'_{\alpha
i})$ where the $p^{m_\alpha}$ is the degree of inseparability of the
field extension $K_{\alpha i(\alpha)}$ over $F$ (see again \cite{W},
p.207). Note that $p^{m_\alpha}$ does not depend on the choice of
the index $i(\alpha)$, because the field $K_{\alpha i(\alpha)}$ is
determined by $\alpha$ up to {\em conjugation} over $F$. From that
point onwards the proof is the same (note in particular that we
shall have $Z'_\alpha=N_{K_{\alpha i(\alpha)}/F}(P'_{
\alpha i},Q'_{\alpha i})$).
\end{Remark}

\medskip

\noindent{\bf Lemma 1.7.2} \, \it For $P,Q\in E(F)$ we have
$\left<Q,P\right>=-\left<P,Q\right>$, i.e., the symbol is
skew-symmetric. \rm

\medskip

\begin{Proof}
Since the symbol is bilinear we have a well-defined homomorphism
$$
\lambda\colon\, E(F)\otimes_{\Z}E(F)\to ST_F(A)\text{ with
}\lambda(P\otimes Q)= \left<P,Q\right>.
$$
Then $\lambda((P+Q)\otimes (P+Q))=(P+Q,P+Q) -(P+Q,0)-(0,P+Q)+(0,0)$.
On the other hand by the bilinearity, it also equals
$\left<P,P\right>+\left<P,Q\right>+
\left<Q,P\right>+\left<Q,Q\right>$. On the diagonal we have
$(P+Q,P+Q)+(0,0)=(P,P)+(Q,Q)$, on $E\times 0$ we have
$(P+Q,0)+(0,0)=(P,0)+(Q,0)$, and on $0\times E$ we have
$(0,P+Q)+(0,0)=(0,P)+(0,Q)$. Putting these facts together we get
$\left<P,Q\right>+\left<Q,P\right>=0$.
\end{Proof}

\medskip

\mysubsection{A useful lemma} \label{subsec1.8}

\medskip

We will often have occasion to use the following simple observation:

\medskip

\noindent{\bf Lemma 1.8.1} \, \it Let $F$ be any field and $\ell$ a
prime. Let $P, Q$ be points in $E(F)$, $Q'\in E(\Fbar)$ s.t. $\ell
Q'=Q$, and put $K=F(Q')$. Consider the statements

\begin{itemize}
\item[{\rm (a)}]
$P\in N_{K/F}E(K)\,\,\mod\,\ell \,E(F)$

\item[{\rm (b)}]
$\left<P,Q\right>\in \ell \,T_F(A)$
\end{itemize}
Then (a) implies (b). \rm

\medskip

\begin{Proof}
Let $P=N_{K/F}(P')+\ell P_1$, with $P'\in E(K)$, $P_1\in E(F)$. Then
$\left<P,Q\right> -\ell \left<P_1,Q\right>
=\left<N_{K/F}(P'),Q\right>$, which equals, by the projection
formula,
$$
N_{K/F}(\left<P',\Res_{K/F}Q\right>) = N_{K/F}(\left<P',\ell
Q'\right>)=\ell N_{K/F} (\left<P',Q'\right>)\in \ell T_F(A).
$$
\end{Proof}

\qed

\bigskip

\section{Symbols, cup products, and $H^2_s(F, E^{\otimes 2})$}
\label{sec2}

\mysubsection{Degeneration of the spectral sequence}
\label{subsec2.1} Notations and assumption are as before. Let $\ell$
be a prime number with $\ell\not=\Char(F)$. We work here with
$\Z_\ell$ coefficients, but the results are also true for
$\Z/\ell^s$ coefficients (any $s\Ge 1$) and
$\Q_\ell$-coefficients.

\noindent{\bf Lemma 2.1.1} \, \it The Hochschild-Serre spectral
sequence
$$
E_2^{pq}=H^p(F,H_{\et}^q(\Abar,\Z_\ell(s)) \Longrightarrow
H_{\et}^{p+q}(A,\Z_\ell(s))
$$
degenerates at $d_2$-level (and at all $d_t$-levels, $t\Ge 2$) for
all $r$. \rm

\medskip

\begin{Remark}
We could take here any abelian variety $A$ instead of $E\times E$.
\end{Remark}

\begin{Proof}
This follows from ``weight'' considerations. Consider on $A$
multiplication by $n$, i.e. $n\colon\, A\to A$ is the map $x\to nx$.
Then we have a commutative diagram
\begin{equation*}
\begin{CD}
H^p(F,H_{\et}^q(\Abar,-))@>{d_2}>>
H^{p+2}(F,H_{\et}^{q-1}(\Abar,-))\\
@V{n^*}VV @V{n^*}VV\\
H^p(F,H_{\et}^q(\Abar,-)) @>{d_2}>>
  H^{p+2}(F,H_{\et}^{q-1}(\Abar,-))
\end{CD}
\end{equation*}
On the left we have multiplication by $n^q$, on the right by
$n^{q-1}$. this being true for any $n>0$, we must have $d_2=0$.
\end{Proof}

\bigskip

\mysubsection{Kummer sequence {\rm (and some notations)}}
\label{subsec2.2}

\medskip

If $E$ is an elliptic curve defined over $F$, we write (by abuse of
notation)
$$
E[\ell^n]:= \ker\{E(\Fbar){\overset{\ell^n}\longrightarrow}
E(\Fbar)\},
$$
and we have the (elliptic) Kummer sequence
\begin{equation*}
\begin{CD}
0@>>> E[\ell^n]@>>> E(\Fbar)@>{\ell^n}>> E(\Fbar)@>>> 0.
\end{CD}
\end{equation*}
This is an exact sequence of $\Gal(\Fbar/F)$-modules and gives us a
short exact sequence of (cohomology) groups
$$
0\longrightarrow E(F)/\ell^n\longrightarrow
H^1(F,E[\ell^n])\longrightarrow H^1(F,\Ebar)[\ell^n]\longrightarrow
0
$$
Taking the limit over $n$, one gets the homomorphism
$$
\delta_\ell^{(1)}\colon\, E(F)\longrightarrow H^1(F,T_\ell(E))
$$
where $T_\ell(E)=\Limleft{n}\,E[\ell^n]$ is the Tate group.

This allows us to define
\begin{align*}
&[\cdot ,\cdot]_\ell\colon\,
  E(F)\otimes E(F)\longrightarrow
H^2(F,T_\ell(E)^{\otimes 2})\\
\intertext{by} &[P,Q]_\ell=\delta_\ell^{(1)}(P)\cup
  \delta_\ell^{(1)}(Q)
\end{align*}
We have similar maps (and we use the same notations) if we take
$E(F)/\ell^n$ and $E[\ell^n]^{\otimes 2}$.

Explicitly the map $\delta_\ell^{(1)}$ is given by the following:
For $P\in E(F)$ the cohomology class $\delta_\ell^1(P)$ is
represented by the $1$-cocycle
$$
\Gal(\Fbar/F) \, \rightarrow \, T_\ell(E),
$$
$$
\sigma \mapsto \left(\sigma\left(
\frac{1}{\ell}\,P\right)-\frac1\ell\,P,  \,
\sigma\left(\frac{1}{\ell^2}\,P\right) -\frac{1}{\ell^2}\,P, \,
\ldots\ldots \, \right).
$$

\bigskip

\mysubsection{Comparison with the usual cycle class map}
\label{subsec2.3}

\medskip

For every smooth projective variety $X$ defined over $F$ there is
the cycle class map to continuous cohomology as defined by Jannsen
(\cite{J}, lemma 6.14)
$$
\cl_\ell^{(i)}\colon\, \CH^i(X)\longrightarrow
H_{\cont}^{2i}(X,\Z_\ell(i))
$$
Taking now $X=E$, resp. $X=A$, and using the degeneration of the
Hochschild-Serre spectral sequence we get
$$
\cl_{\ell}^{(1)}\colon\, CH_{(0)}^1(E) \longrightarrow
H^1(F,H_{\et}^1(\Ebar, \Z_\ell(1)))
$$
where $\CH_{(0)}^1(E)$ is the Chow group of $0$-cycles on $E$ of
degree $0$, resp.
$$
\cl_\ell^{(2)}\colon\,T_F(A)\longrightarrow
H^2(F,H_{\et}^2(\Abar,\Z_\ell(2))).
$$

\noindent{\bf Lemma 2.3.1} \, \it There is a commutative diagram
\begin{equation*}
\begin{CD}
E(F) @>{\delta_\ell^{(1)}}>>
H^1(F,T_\ell(E))\\
@VVV @VV{\cong}V\\
\CH_{(0)}^1(E) @>{cl_\ell^{(1)}}>>
H^1(F,H_{\et}^1(\Ebar,\Z_\ell(1)))
\end{CD}
\end{equation*}
where the vertical map on the left is $P\mapsto (P)-(0)$ and the one
on the right comes from the well-known isomorphism
$T_\ell(E)\arrowsim H_{\et}^1(\Ebar,\Z_\ell(1))$. \rm

\medskip

\begin{Proof}
See \cite{R}, proof of the lemma in the appendix.
\end{Proof}

\bigskip

\mysubsection{The symbolic part of cohomology} \label{subsec2.4}

\medskip

Let $K/F$ be any finite extension. Given a pair of points $P, Q \in
E(K)/\ell$, with associated classes $\delta_\ell^{(1)}(P),
\delta_\ell^{(1)}(Q)$ in $H^1(K,E[\ell])$. We have seen in section
2.2 that by taking cup product in Galois cohomology, we get a class
$[P,Q]_\ell$ in $H^2(K, E[\ell]^{\otimes 2})$. By taking the norm
(corestriction) from $H^2(K,E[\ell]^{\otimes 2})$ to
$H^2(F,E[\ell]^{\otimes 2})$, we then get a class
$$
N_{K/F}([P,Q]_\ell) \, \in H^2(F, E[\ell]^{\otimes 2}).
$$

We define the {\it symbolic part} $H^2_s(F,E[\ell]^{\otimes 2})$ of $H^2(F,E[\ell]^{\otimes 2})$ to
be the $\F_\ell$-subspace generated by such {\it norms of symbols}
$N_{K/F}([P,Q]_\ell)$, where $K$ runs over all possible finite
extensions of $F$ and $P, Q$ run over all pairs of points in
$E(K)/\ell$.

Note the similarity of this definition with the description of
$T_F(E \times E)/\ell$ via Lemma 1.7.1.

\bigskip

\mysubsection{Summary} \label{subsec2.5}

\medskip

The remarks and maps of the previous sections can be subsumed in the
following:

\medskip

\noindent{\bf Proposition 2.5.1} \, \it There exist maps and a
commutative diagram
\begin{equation*}
\left<\,\,,\,\,\right>
\lower77pt\hbox{\epsfig{file=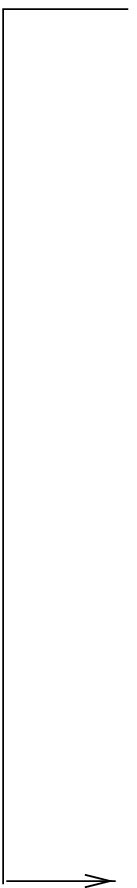,height=6.9 true cm}}
\begin{CD}
E(F)^{\otimes 2} @>{\delta_\ell^{(1)}\otimes
  \delta_\ell^{(1)}}>>
H^1(F,T_\ell(E))^{\otimes 2}\\
@VVV @VV{\cong}V\\
\CH_{(0)}^1(E)^{\otimes 2} @>{\cl_\ell^{(1)}
  \otimes \cl_\ell^{(1)}}>>
H^1(F,H_{\et}^1(\Ebar,\Z_\ell(1)))^{\otimes 2}\\
@VV{p_1^*\otimes p_2^*}V \raise5pt\hbox{${{\displaystyle\hookdownarrow\atop
\displaystyle H_{\cont}^2
  (E,\Z_\ell(1))^{\otimes 2}}\atop
{\big\downarrow \rlap{$\vcenter{\hbox{$\scriptstyle
  p_1^*\otimes p_2^*$}}$}}}$}\\
\CH^1(A)^{\otimes 2} @>{\cl_\ell^{(1)} \otimes\cl_\ell^{(1)}}>>
H_{\cont}^2(A,\Z_\ell(1))^{\otimes 2}\\
@VV{\cap}V @VV{\cup}V\\
\CH^2(A) @>{\cl_\ell^{(2)}}>>
H_{\cont}^4(A,\Z_\ell (2))\\
\hookuparrow @. \hookuparrow\\
T_F(A) @>{\cl_\ell^{(2)}}>> H^2(F,H_{\et}^2
  (\Abar,\Z_\ell(2)))\\
&&\hookuparrow\\
&\vbox{\beginpicture \setcoordinatesystem units <.50cm,.50cm>
\linethickness=1pt \setshadesymbol ({\thinlinefont .}) \setlinear
%
%
\put{\lower1pt\hbox{$c$}} [lB] at  9.811 17.486
%
%
\put{\lower2pt\hbox{$\scriptstyle\ell$}} [lB] at 10.128 17.264
%
%
\linethickness= 0.500pt \setplotsymbol ({\thinlinefont .}) \plot
8.255 18.812 12.922 16.932 /
%
%
\plot 12.660 16.928 12.922 16.932 12.689 17.181 /
\linethickness=0pt \putrectangle corners at  8.230 17.837 and 12.948
16.701
\endpicture}
&H^2(F,H_{\et}^1(\Ebar,\Z_\ell(1))^{\otimes 2})
\end{CD}
\lower122pt \hbox{\epsfig{file=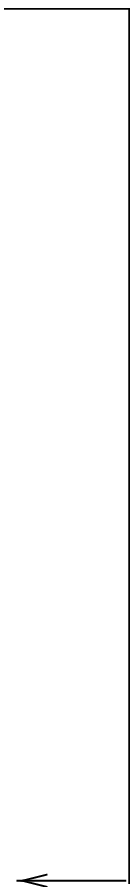,height=6.90 true cm}}
\,\,\lower10pt\hbox{$\kern4pt$}
\end{equation*}
where $c_\ell$ is defined via the K\"unneth formula
$H_{\et}^2(\Abar,\Z_\ell(2))= H_{\et}^2(\Ebar,\Z_\ell(2))\oplus
H_{\et}^1(\Ebar,\Z_\ell(1))^{\otimes 2} \oplus
H_{\et}^2(\Ebar,\Z_\ell(2))$ and the injection on the relevant
($=$ middle) term. In particular
$[P,Q]_\ell:=\delta_\ell^{(1)}(P)\cup
\delta_\ell^{(1)}(Q)=c_\ell(\left<P,Q\right>)$ for $P,Q\in E(F)$.
Moreover the same holds for $\Z/\ell^n$-coefficients (instead of
$\Z_\ell$) and also for $\Q_\ell$-coefficients. \rm

\medskip

\noindent{\bf Remark}: \, Note that the map $\langle , \rangle$ on
the left from $E(F)^{\otimes 2}$ to $T_F(A)$ is the symbol map, and
that the map from $H^1(F, H^1_{\rm et}(\overline E,
\Z_\ell(1)))^{\otimes 2}$ to $H^2(F, H^1_{\rm et}(\overline E,
\Z_\ell(1))^{\otimes 2})$ on the right is the one given by taking the cup product
in group cohomology.

\medskip

\noindent{\bf Proof} (and further explanation of the maps)

The map $p_1^*\otimes p_2^*$ is induced by the projections
$p_i\colon\,A=E\times E\to E$. The commutativity in upper rectangle
comes from Lemma 2.3.1. The other rectangles are all ``natural''.
\qed

\medskip

\noindent{\bf Corollary 2.5.2} \, \it With respect to the
decomposition
$$
H_{\et}^1(\Ebar,\Z_\ell(1))^{\otimes 2}= T_\ell(E)^{\otimes
2}\simeq\Sym^2T_\ell(E) \oplus\Lambda^2 T_\ell(E)
$$
we have that $[P,Q]_\ell=c_\ell(\left<P,Q\right>) \in
H^2(F,\Sym^2T_\ell(E))$ if $\ell\not=2$ (and similarly for
$\Z/\ell^s$-coefficients). \rm

\medskip

\noindent{\bf Proof}

\noindent {\bf Step 1.}\enspace For the sake of clarity of the proof
we shall first take two different elliptic curves $E_1$ and $E_2$
(or, if one prefers, two ``different copies'' $E_1$ and $E_2$ of
$E$). Put $A_{12}=E_1\times E_2$. There exist obvious analogs of the
maps and the commutation diagram of Proposition 2.5.1; only -- of
course -- one should now write $E_1(F)\times E_2(F)$, etc. In
particular we have again the cup product on the right:
$$
H^1(F,T_\ell(E_1))\otimes H^1(F,T_\ell(E_2))
{\overset\cup\longrightarrow}H^2(F,T_\ell(E_1) \otimes T_\ell(E_2))
$$
where we write $T_\ell(E_1)=H_{\et}^1(\Ebar_i,\Z_\ell(1)),\,
i=1,2$. We get $c_\ell(\left<P,Q\right>)=\delta_\ell^{(1)}
(P)\cup\delta_\ell^{(1)}(Q)$. Now consider also $A_{21}=E_2\times
E_1$ and the corresponding diagram for $A_{21}$. Consider the
natural isomorphism $t\colon\,A_{12}\to A_{21}$ given by
$t(x,y)=(y,x)$ for $x\in E_1$ and $y\in E_2$ and also the
corresponding isomorphism
$$
t_*\colon\, H^1(\Ebar_1)\otimes H^1(\Ebar_2) \longrightarrow
H^1(\Ebar_2)\otimes H^1(\Ebar_1).
$$

\bigskip\noindent
{\bf Claim}.\enspace \it If $P\in E_1(F)$ and $Q\in E_2(F)$ then
$$
c_{\ell,A_{21}}(\left<Q,P\right>)= -t_*(\delta_\ell^{(1)}(P)\cup
\delta_\ell^{(1)}(Q))
$$
\rm

\medskip
Here we have written -- in order to avoid confusion --
$c_{\ell,A_{21}}$ for $c_\ell$ in the diagram relative to $A_{21}$.

\bigskip

\noindent Proof of the claim:\enspace
$$
c_{\ell,A_{21}}(\left<Q,P\right>)=
\delta_\ell^{(1)}(Q)\cup\delta_\ell^{(1)}(P)=
-t_*\left(\delta_\ell^{(1)}(P)\cup \delta_\ell^{(1)}(Q)\right)
$$
where the first equality is from the diagram (for $A_{21}$) in
Proposition 2.5.1 and the second equality is a well-known property
in cohomology (see for instance \cite{Brown}, p. 111, (3.6), or
\cite{Spanier}, chap.5, \S6, p. 250).

\bigskip\noindent
{\bf Step 2.}\enspace Returning to the case $E_1=E_2=E$, we have
$\left<Q,P\right>=-\left<P,Q\right>$ in $T_F(A)$ by Lemma 1.7.2.

\bigskip\noindent
{\bf Step 3.}\enspace Let us write $c_\ell(\left<P,Q\right>)=\alpha+\beta$,
with $\alpha$ in $H^2(F,\Sym^2T_\ell(E))$ and $\beta$ in
$H^2(F,\Lambda^2T_\ell(E))$. We get $c_\ell(\left<P,Q\right>)=-c_\ell
\left<Q,P\right>)$ from Step 2, and next from Claim in Step 1 we
have $-c_\ell(\left<Q,P\right>)=t_*(\delta_\ell^{(1)}
(P)\cup\delta_\ell^{(1)}(Q))$, hence
$\alpha+\beta=t_*(\alpha+\beta)=\alpha-\beta$, hence $\beta=0$ if
$\ell\not=2$.

\qed

\bigskip

\noindent{\bf Definition 2.5.3} \, \it Let
$$
s_\ell: \, T_F(E \times E)/\ell \, \rightarrow \, H^2(F,
E[\ell]^{\otimes 2})
$$
be the homomorphism induced by the reduction of $c_\ell$ modulo
$\ell$ and by using the isomorphism of $H^1(E_{\overline F},
\Z_\ell(1))$ with the $\ell$-adic Tate module of $E$, which is the
inverse limit of $E[\ell^n]$ over $n$. \rm

\medskip

Thanks to the discussion above, as well as the definition of the symbolic part
of $H^2(F, E[\ell]^{\otimes 2})$ in section 2.4, and Lemma 1.7.1, we
obtain the following:

\medskip

\noindent{\bf Proposition 2.5.4} \, \it The image of $T_F(E \times
E)/\ell$ under $s_\ell$ is $H^2_s(F, E[\ell]^{\otimes 2})$. \rm

\bigskip

\section{A key Proposition}
\label{sec3}

\medskip

Let $E/F$ be an elliptic curve over a local field $F$ as in Theorem A,
i.e., non archimedean and with finite
residue field of characteristic $p>2$, and with
$E$ having good, ordinary
reduction. We will henceforth
take $\ell=p$.

A basic fact is that the representation $\rho_p$
of $G_F$ on $E[p]$ is reducible, so the
matrix of this representation is triangular. Since the
determinant is the mod $p$ cyclotomic character $\chi_p$, we may write
$$
\rho_p \, = \, \begin{pmatrix}\chi_p\nu^{-1} & \ast\\ 0 &
\nu\end{pmatrix},\leqno(3.1)
$$
where $\nu$ is an unramified character of finite order, such that
$E[p]$ is semisimple (as a $G_F$-module) iff $\ast \, = \,
0$. Note that $\nu$ is necessarily of order at most $2$ when $E$ has
multiplicative reduction, with $\nu=1$ iff $E$ has {\it split}
multiplicative reduction. On the other hand, $\nu$ can have
arbitrary order (dividing $p-1$) when $E$ has good, ordinary
reduction. In any case, there is a natural $G_F$-submodule $C_F$
of $E[p]$ of dimension $1$, such that we have a short exact
sequence of $G_F$-modules:
$$
0 \to C_F \to E[p] \to C_F^\prime \to 0,\leqno(3.2)
$$
with $G_F$ acting on $C_F$ by $\chi_p\nu^{-1}$ and on
$C_F^\prime$ by $\nu$. Clearly, $E[p]$ is semisimple iff the
sequence $(3.2)$ splits.

\medskip

The natural $G_F$-map $C_F^{\otimes 2} \to E[p]^{\otimes 2}$
induces a homomorphism
$$
\gamma_F: \, H^2(F, C_F^{\otimes 2}) \, \rightarrow \, H^2F,
E[p]^{\otimes 2}).\leqno(3.3)
$$
The key result we prove in this section is the following:

\medskip

\noindent{\bf Proposition C} \, \it Let $F$ be a non-archimedean
local field with odd residual characteristic $p$, and
$E$ an elliptic curve over $F$ with good,
ordinary reduction.
Denote by $Im(s_p)$ the image of $T_F(E \times E)/p$ under the
Galois symbol map $s_p$ into $H^2F, E[p]^{\otimes 2})$. Then
we have
\begin{enumerate}
\item[(a)]
$$
{\rm Im}(s_p) \, \subset \, {\rm Im}(\gamma_F).
$$
\item[(b)]The dimension of ${\rm Im}(s_p)$ is at most $1$, and it
is zero dimensional if either $\mu_p \not\subset F$ or $\nu^2\ne
1$.
\end{enumerate}\rm

\medskip

\noindent{\bf Remark}: \, If $E/F$ has multiplicative reduction, with $\ell$ an
arbitrary odd prime (possibly equal to $p$),
then again the Galois representation $\rho_\ell$ on $E[\ell]$ has a similar shape,
and in fact, $\nu$ is at most quadratic, reflecting the fact that $E$ attains split
multiplicative reduction over at least a quadratic extension of $F$, over which the
two tangent directions at the node are rational. An analogue of Proposition C holds
in that case, thanks to a key result of W.McCallum ([Mac], Prop.3.1), once we assume that $\ell$
does not divide the order of the component group of the special fibre
$\ScrE_{s}$ of the N\'eron model $\ScrE$. For the sake of brevity, we
are not treating this case here.

\medskip

The bulk of this section will be involved in proving the following
result, which at first seems weaker than Proposition C:

\medskip

\noindent{\bf Proposition 3.4} \, \it Let $F, E, p$ be as in
Proposition C. Then, for all points $P, Q$ in $E(F)/p$, we have
$$
s_p(\langle P, Q\rangle) \, \in \, {\rm Im}(\gamma_F).
$$
\rm

\medskip

\noindent{\bf Claim 3.5} \, {\it Proposition 3.4 $\implies$ Part (a) of Proposition C}:

\medskip

Proof of this claim goes via some lemmas.

\medskip

\noindent{\bf Lemma 3.6} \, (Behavior of $\gamma$ under finite extensions)
\, \it Let $K/F$ be finite. We then have two commutative diagrams, one for
the norm map $N=N_{K/F}$ and the other for the restriction map Res$=$Res$_{K/F}$:
\begin{equation*}
\begin{matrix}
H^2(K, C_K^{\otimes 2}) &\longmapright{25}{\gamma_K}
  &H^2(K,E[p]^{\otimes 2})\\
\lmapdown{N}\uparrow{Res}
&&\lmapdown{N}\uparrow{Res}\\
H^2(F, C_F^{\otimes 2}) &\longmaprightsub{25}{\gamma_F}
&H^2(F,E[p]^{\otimes 2})
\end{matrix}
\end{equation*}
\rm

\medskip

This Lemma follows from the compatibility of the exact sequence (3.2) (of Galois modules)
with base extension.

\medskip

\noindent{\bf Lemma 3.7} \, \it In order to prove that Im$(s_p) \subset {\rm Im}(\gamma_F)$, it suffices to
prove it for the image of symbols, i.e., that
$$
{\rm Im}\left(s_p(ST_{F,p}(A))\right) \, \subset \, {\rm Im}(\gamma_F),
$$
where
$$
ST_{F,p}(A): = \, {\rm Im}\left(ST_{F}(A)\rightarrow T_F(A)/p\right).
$$
\rm

\medskip

{\it Proof}. \, Use the commutativity of the diagram in Lemma 3.6 for the norm. \qed

\medskip

This proves {\it Claim 3.5}.

\qed

\medskip

\noindent{\bf Lemma 3.8} \, \it In order to prove that Im$(s_p)
\subset {\rm Im}(\gamma_F)$, we may assume that $\mu_p \subset F$
and that $\nu=1$, i.e., that we have the following exact sequence
for groupschemes over $S$:
$$
0 \to \mu_{p, S} \to \ScrE[p] \to (\Z/p)_S \to 0.\leqno(\ast)
$$
\rm

\medskip

{\it Proof}. \, There is a finite extension $K/F$ such that $(\ast)$ holds over $K$, with $p \nmid [K:F]$.
Now use the diagram(s) in Lemma 3.6 as follows: Let $P, Q \in E(F)\subset E(K)$, then Res$(P)=P$,
Res$(Q)=Q$, and we have
$$
[K:F]s_{F,p}(\langle P, Q\rangle)  =  N\{{\rm Res}(s_{F,p}(\langle P, Q\rangle\},
$$
which equals
$$
N\{s_{K,p}\left(\langle{\rm Res}(P), {\rm Res}(Q)\rangle\right)\} = N\{s_{K,p}(\langle P, Q\rangle)\}.
$$
Therefore, if $s_{K,p}(\langle P, Q\rangle) \subset {\rm Im}(\gamma_K)$, then (again by
using Lemma 3.6 for norm) we see that
$N\{s_{K,p}(\langle P, Q\rangle)\}$ is contained in ${\rm Im}(\gamma_F)$. Hence $[K:F]s_{F,p}(\langle P, Q\rangle)$
lies in Im$(\gamma_F)$. Finally, since $p\nmid [K:F]$, $s_{F,p}(\langle P, Q\rangle)$
itself belongs to Im$(\gamma_F)$.

\qed

\medskip

\subsection*{$3.9$ \, Proof of Proposition 3.4}

\medskip

Since we may take $\mu_p\subset F$ and $\nu=1$, the exact sequence
$(\ast)$ in Lemma 3.8 holds, compatibly with the
corresponding one over $F$ of Gal$(\overline F/F)$-modules. Taking
cohomology, we get a commutative diagram of $\F_p$-vector spaces
with exact rows:
\begin{equation*}
\begin{matrix}
\ScrO_F^\ast/p &\to
&{E(F)/p}&\to &{\Z/p}\\
\downarrow&
&\downarrow &&\downarrow\\
\overline H^1(F,\mu_p) &\to &{\overline H^1(F,E[p])} &\to
&{\overline H^1(F,\Z/p)} &\to
&{0}
\end{matrix}
\end{equation*}

Here the vertical maps are isomorphisms (see the Appendix), which induce the horizontal maps on the top row.

\medskip

\noindent{\bf Definition 3.10} \, \it Put
$$
U_F:= \, {\rm Im}\left(\ScrO_F^\ast/p \, \rightarrow \, E(F)/p\right),
$$
and choose a non-canonical decomposition
$$
E(F)/p \, \simeq \, U_F \oplus W_F, \quad {\rm with} \quad W_F \simeq \Z/p.
$$
\rm

\medskip

\noindent{\bf Notation 3.11} \, If $S_1, S_2$ are subsets of $E(F)/p$, we denotes by
$\langle S_1, S_2\rangle$ the subgroup of $ST_{F,p}(A)$ generated by the symbols
$\langle s_1, s_2\rangle$, with $s_1\in S_1$ and $s_2\in S_2$.

\medskip

\noindent{\bf Lemma 3.12} \, \it $ST_{F,p}(A)$ is generated by the two vector subspaces
$$
\Sigma_1:= \, \langle U_F, E(F)/p\rangle \quad {\rm and} \quad
\Sigma_2:= \, \langle E(F)/p, U_F\rangle.
$$
\rm

\medskip

{\it Proof}. \, $ST_{F,p}(A)$ is clearly generated by $\Sigma_1, \Sigma_2$ and by
$\langle W_F, W_F\rangle$. However, $W_F$ is one-dimensional and the pairing
$\langle \cdot, \cdot\rangle$ is skew-symmetric, so $\langle W_F, W_F\rangle = 0$.

\qed

\medskip

Now we have a commutative diagram

\noindent(3.13-i)

\begin{equation*}
\begin{matrix}
\ScrO_F^\ast/p \otimes E(F)/p &\mapright{\alpha_1}
&E(F)/p \otimes E(F)/p\\
\rmapdown{\sigma_1} & & \rmapdown{s_p}\\
H^2(F, \mu_p\otimes E[p]) &\mapright{\beta_1} &H^2(F, E[p]^{\otimes 2}),
\end{matrix}
\end{equation*}
where the top map $\alpha_1$ factors as
$$
\ScrO_F^\ast/p \otimes E(F)/p \rightarrow U_F\otimes E(F)/p \rightarrow E(F)/p \otimes E(F)/p.
$$

We also get a similar diagram (3.13-ii) by replacing $\ScrO_F^\ast/p \otimes (E(F)/p)$
(resp. $U_F\otimes (E(F)/p)$) by $(E(F)/p) \otimes \ScrO_F^\ast/p$ (resp.
$(E(F)/p) \otimes U_F$). The maps $\alpha_j$ and their factoring are obvious, $s_p$ is the map
constructed in section 2, and the vertical maps $\sigma_j$ are defined entirely analogously.

\medskip

\noindent{\bf Lemma 3.14} \, \it The image of $s_p: ST_{F,p}(A) \rightarrow H^2(F, E[p]^{\otimes 2})$
is generated by $\beta_1(Im (\sigma_1))$ and $\beta_2(Im (\sigma_2))$. \rm

\medskip

{\it Proof}. \, Immediate by Lemma 3.12 together with the commutative diagrams (3.13-i) and (3.13-ii).

\qed

\medskip

Tensoring the exact sequence
$$
0\to\mu_p\to E[p]\to\Z/p\to 0\leqno(3.15)
$$
with $\mu_p$ from the left and the right, and taking Galois cohomology, we get two natural
homomorphisms
$$
\gamma_1: H^2(F, \mu_p^{\otimes 2}) \, \rightarrow \, H^2(F, \mu_p \otimes E[p])\leqno(3.16-i)
$$
and
$$
\gamma_2: H^2(F, \mu_p^{\otimes 2}) \, \rightarrow \, H^2(F, E[p] \otimes \mu_p).\leqno(3.16-ii)
$$

\medskip

\noindent{\bf Lemma 3.17} \, \it ${\rm Im}(\sigma_j) \, \subset \, {\rm Im}(\gamma_j)$, for $j=1, 2$. \rm

\medskip

{\it Proof of Lemma 3.17}. \, We give a proof for $j=1$ and leave the other (entirely similar) case to the reader.
By tensoring (3.15) by $\mu_p$, we obtain the following exact sequence of $G_F$-modules:
$$
0\to\mu_p^{\otimes 2}\to \mu_p\otimes E[p]\to\mu_p\to 0\leqno(3.18)
$$

Consider now the following commutative diagram, in which the bottom row is exact:

\noindent(3.19)

\begin{equation*}
\begin{matrix}
&& U_F\otimes E(F)/p &&\\
& & \rmapdown{\tilde\sigma_1} & &\\
&&{\overline H}^2(F, \mu_p\otimes E[p]) &\mapright{\varepsilon_0} &{\overline H}^2(F, \mu_p)\\
&&\rmapdown{i_1} & &\rmapdown{} &&\\
H^2(F, \mu_p^{\otimes 2}) &\mapright{\gamma_1} &H^2(F, \mu_p\otimes E[p]) &\mapright{\varepsilon}
&H^2(F, \mu_p)
\end{matrix}
\end{equation*}

\noindent{}where $\sigma_1=i_1\circ {\tilde \sigma}_1$ is the map from (3.13-i).

\medskip

To begin, the exactness of the bottom row follows immediately from the exact sequence (3.18).
The factorization of $\sigma_1$ is the crucial point, and this holds because $U_F$ comes from $\ScrO_F^\ast/p$
and is mapped to $H^1(F, \mu_p)$ via $\Hbar^1(F, \mu_p)$, which is the image of $H^1_{\rm fl}(S, \mu_p)$ (see Lemma A.3.1 in the Appendix). Similarly, $E(F)/p \simeq \ScrE(S)/p$ maps into $\Hbar^1(F, E[p])$, the image of $H^1_{\rm fl}(S, \ScrE[p])$; therefore the tensor product maps into $\Hbar^2(F, \mu_p\otimes E[p])$, the image of $H^2_{\rm fl}(S, \mu_p\otimes \ScrE[p])$.

To prove Lemma 3.17, it suffices, by the exactness of the bottom row, to see that Im$(\sigma_1)$ is contained
in Ker$(\varepsilon)$, i.e., to see that $\varepsilon\circ \sigma_1 = 0$. Luckily for us, this
composite map factors through $\varepsilon_0\circ {\tilde \sigma_1}$, which vanishes because $H^2_{\rm fl}(S, \mu_{p, S})$,
and hence $\Hbar^2(F, \mu_p)$, is zero by \cite{Mi2}, part III, Lemma 1.1.

\qed

\medskip

Putting these Lemmas together, we get the truth of Proposition 3.4.

\medskip

\subsection*{${\bf 3.20}$ \, Proof of Proposition C}

\medskip

As we saw earlier, Proposition 3.4, which has now been proved,
implies (by Claim 3.5) part (a) of Proposition C. So we need to prove only part
(b).

By part (a) of Prop. C, the dimension of the image of $s_p$ is at
most that of $H^2(F, C^{\otimes 2})$. Since $E[p]$ is selfdual, the
short exact sequence (3.2) shows that the Cartier dual of $C$ is
$C'$. It follows easily that $(C^{\otimes 2})^D$ is ${C'}^{\otimes
2}(-1)$. By the local duality, we then get
$$
H^2(F, C^{\otimes 2}) \, \simeq \, H^0(F, {C'}^{\otimes 2}(-1))^\vee,
$$
As $C'$ is a line over $\F_p$ with $G_F$-action, the dimension of the group on the right is less than or equal to $1$,
with equality holding iff ${C'}^{\otimes 2} \simeq \mu_p$.

Since $G_F$ acts on $C'$ by an unramified character $\nu$, for ${C'}^{\otimes 2}$ to be $\mu_p$ as a $G_F$-module,
it is necessary that $\mu_p\subset F$. So
$$
\mu_p \not\subset F \, \implies \, {\rm Im}(s_p)=0.
$$

Now suppose $\mu_p \subset F$. Then for Im$(s_p)$ to be non-zero, it is necessary that ${C'}^{\otimes 2} \simeq \Z/p$, implying that
$\nu^2=1$.

\qed

\bigskip

\section{Vanishing of $s_p$ in the wild case}
\label{sec4}

\medskip

We start with a simple Lemma, doubtless known to experts.

\medskip

\noindent{\bf Lemma 4.0} \, \it Let $E/F$ be ordinary. Then the Galois representation on $E[p]$ is one of the following (disjoint) types:
\begin{enumerate}
\item[(i)]The  wild inertia group $I_p$ acts non-trivially, making $E[p]$ non-semisimple;
\item[(ii)]$E[p]$ is semisimple, in which case $p \nmid F(E[p]):F]$;
\item[(iii)]$E[p]$ is unramified and non-semisimple, in which case $[F(E[p]):F]=p$.
\end{enumerate}
\rm

\medskip

In case (i), we will say that we are in the {\it wild case}.

\medskip

{\it Proof}. \, Clearly, when $I_p$ acts non-trivially, the triangular nature of the representation on $E[p]$ makes $E[p]$ non-semisimple. Suppose from here on that $I_p$ acts trivially. If $E[p]$ is semisimple, $F(E[p])/F$ is necessarily a prime-to-$p$ extension, again because the representation is solvable in the ordinary case. This leaves the final possibility when $E[p]$ is non-semisimple, but $I_p$ acts trivially. We claim that the tame inertia group $I_t$ also acts trivially. Since $I_t \simeq \lim_m \F_{p^m}^\ast$, its image in GL$_2(\F_p)$ ($\simeq$ Aut$(E[p])$) will have prime-to-$p$ order, and so must be (conjugate to) a subgroup $H$ of the diagonal group $D \simeq (\F_p^\ast)^2$. $H$ cannot be central, as the semisimplification is a direct sum of two characters, at most one of which is ramified, and for the same reason, $H$ cannot be all of $D$. It follows that, when $I_t$ acts non-trivially, $H$ must be conjugate to
$$
\left\{\begin{pmatrix}a&0\\0&1\end{pmatrix}\, \vert \, a \in C \subset \F_p^\ast\right\},
$$
for a subgroup $C \ne \{1\}$ of $\F_p^\ast$.
An explicit calculation shows that any element of GL$_2(\F_p)$ normalizing $H$ must be contained in the group generated by $D$ and the Weyl element $w=\begin{pmatrix}0 & -1\\ 1 & 0\end{pmatrix}$. Again, since the representation is triangular, $w$ cannot be in the image. Consequently, if $I_t$ acts non-trivially, the image of $\ScrG_k$ must also be in $D$, implying that $E[p]$ is semisimple, whence the claim. Thus in the third (and last) case, $E[p]$ must be unramified and non-semisimple. As $\ScrG_k$ is pro-cyclic, the only way this can happen is for the image to be unipotent of order $p$.
\qed

\medskip

\noindent{\bf Remark}: \, An example of Rubin
shows that for a suitable, ordinary elliptic curve $E$ over a $p$-adic field,
we can be in case (iii). To this,
start with an example $E/F$ in case (i), with E[p] unipotent and non-trivial,
which is easy to produce.  Let $L$ be the ramified extension  of $F$
of degree p obtained by trivializing $E[p]$. Let $K$ be the unique unramified extension
of $F$ of degree $p$, and let $M$ be an extension of $F$ of degree $p$ inside
the compositum $LK$, different from both $L$ and $K$.  Then over $M$, $E[p]$
has the form we want.  Since M does not
contain $L$, $E[p]$ is a non-trivial $\ScrG_M$-module. On the other hand, this representation
(of $\ScrG_M$) is unramified  because all the $p$-torsion is defined over the
unramified extension $MK$ of $M$, which contains $L$.

\medskip

The object of this section is to prove the following:

\medskip

\noindent{\bf Proposition D} \, \it Suppose $E/F$ is an elliptic
curve (with good ordinary reduction) over a non-archimedean local
field $F$ of odd residual characteristic $p$. Assume that we are in the wild case. Then we have
$$
s_p(T_F(E \times E)/p) \, = \, 0.
$$
\rm

\medskip

Combining this with Proposition C, we get the following

\medskip

\noindent{\bf Corollary 4.1} \, \it Let $F$ be a non-archimedean
local field with residual characteristic $p > 2$, and
$E$ an elliptic curve over $F$ with good,
ordinary reduction. If Im$(s_p)$ is non-zero, then either
$[F(E[p]):F] \leq 2$ or $F(E[p])$ is an unramified $p$-extension. \rm

\medskip

Indeed, if Im$(s_p)$ is non-zero,  Proposition D says that $I_p$ acts trivially on
$E[p]$. If the Gal$(\overline F/F)$- action on $E[p]$ is semisimple, then by Proposition C,
$[F(E[p]):F] \leq 2$. Since the tame inertia group has order prime to $p$,
the only other possibility, thanks to Lemma 4.0, is for $E[p]$ to be non-semisimple and unramified. As the residual Galois group
is cyclic, this also forces $[F(E[p]):F]=p$. Hence the Corollary (assuming the truth of Proposition D).

\medskip

{\it Proof of Proposition D}. \, Let us first note a few basic
things concerning base change to a finite extension $K/F$ of degree
$m$ prime to $p$. To begin, since $E/F$ has ordinary reduction, the
Galois representation $\rho_F$ on $E[p]$ is triangular, and it is
semisimple iff the image does not contain any element of order $p$.
It follows that $\rho_K$ is semisimple iff $\rho_F$ is semisimple.
Moreover, the functoriality of the Galois symbol map relative to the
respective norm and restriction homomorphisms, together with the
fact that the composition of restriction with norm is multiplication
by $m$, implies, as $p\nmid m$, that for any $\theta\in T_F(E\times
E)/p$, we have
$$
s_{p,K}({\rm res}_{K/F}(\theta)) = 0 \, \implies \, s_p(\theta)=0.
$$

So it suffices to prove Proposition D, after possibly replacing $F$
by a finite prime-to-$p$ extension, under the assumption that $F$
contains $\mu_p$ {\it and} $\nu=1$, still with $E[p]$
non-semisimple and wild. Thus we have a non-split, short exact sequence of
finite flat groupschemes over $S$:
$$
0\to\mu_{p,S}\to\ScrE[p]\to (\Z/p)_S\to 0,\leqno(4.2)
$$
with the representation $\rho_F$ of $G = $Gal$(\overline F/F)$
on $E[p]$ having the form:

\noindent{(4.3)}
$$
\begin{pmatrix}
1 &\alpha\\
0 &1
\end{pmatrix}
$$
relative to a suitable basis. Here $\alpha: G \rightarrow \Z/p$
is a non-zero homomorphism, and $F(\alpha)$, the smallest extension
of $F$ over which $\alpha$ becomes trivial, is a ramified
$p$-extension.

Using the exact sequence (4.2), both over $S$ and over $F$, we get
the following commutative diagram:

\noindent{(4.4)}
\begin{equation*}
\begin{matrix}
\ScrE[p](S) &\mapright{0} &\Z/p &\to &H_{\fl}^1(S,\mu_p)
&\mapright{\psi_S}
&H_{\fl}^1(S,\ScrE[p])&\twoheadrightarrow &H^1_{\fl}(S, \Z/p) \\
\vspace{5pt} \Big\Vert &&\Big\Vert
&&\Big\downarrow &&\Big\downarrow & &\Big\downarrow \\
\vspace{5pt} E[p](F) &\mapright{0} &\Z/p &\to &H^1(F,\mu_p)
&\mapright{\psi} &H^1(F,E[p]) &\rightarrow &H^1(F,\Z/p)
\\
\vspace{5pt}
&&&&\Big\Vert &&\cap&&\\
&&&&X+Y &&{Z = \psi(Y)}&&
\end{matrix}
\end{equation*}
where $X\subset\Hbar^1(F,\mu_p)\subset H^1(F,\mu_p)$ is the subspace
given by the image of $\Z/p$. Since $\mu_p$ is in $F$, we can
identify $H^2(F,\mu_p^{\otimes 2})$ with $\Br_F[p]\simeq \Z/p$. By
our assumption, $X \ne 0$ as the sequence (4.2) does {\em not}
split. Now take $e\in X$, with $e\not=0$. Then
$$
e \, \in \ScrO_F^\ast/p \, \subset \, F^\ast/{{F^\ast}^p} \, = \,
H^1(F, \mu_p).
$$
Since $K = F(e^{1/p})$ is clearly the smallest extension of $F$ over
which (4.2) splits, $K$ is also $F(\alpha)$ and hence a ramified
$p$-extension of $F$, since the Galois representation on $E[p]$ is wild.
Then, applying \cite{Serre3}, Prop.5 (iii), we
get a $v\in\Hbar^1(F,\mu_p)$ which is not a norm from $K$, and so
the cup product $\{e,v\}$ is non-zero in $\Br_F[p]$ (cf.
\cite{Tate3}, Prop.4.3).

Now consider the commutative diagram \noindent{(4.5)}
\begin{equation*}
\begin{CD}
H^1(F,\mu_p)^{\otimes 2} @>{\cup}>>
  H^2(F,\mu_p^{\otimes 2})\\
@V{\psi^{\otimes 2}}VV @VV{\gamma}V\\
H^1(F,E[p])^{\otimes 2}
  @>{\cup}>> H^2(F,E[p]^{\otimes 2})
\end{CD}
\end{equation*}
where $\psi$ is the map defined in the previous diagram, and
$\gamma=\gamma_F$ is the map defined in section 3 with $C_F=\mu_p$
and $C_F'=\Z/p$.  By Proposition C, the image of $s_p$ is contained
in that of $\gamma$. On the other hand, $\psi(e) \cup \psi(v) = 0$
because $\psi(e) = 0$. Hence $\gamma(\{e,v\}) = 0$, and the image of
$s_p$ is zero as asserted.

\bigskip

\section{Non-trivial classes in Im$(s_p)$ in the non-wild case with $[F(E[p]):F]'\leq 2$}
\label{sec5}

\medskip

For any finite extension $K/F$, we will write $[K:F]'$ for the prime-to-$p$ degree of $K/F$.

\medskip

\noindent{\bf Proposition E} \, \it Let $E$ be an elliptic curve
over a non-archimedean local field $F$ of residual characteristic $p
\ne 2$. Assume that $E$ has good ordinary reduction, with trivial action on $I_p$ on $E[p]$,
such that
\begin{enumerate}
\item[(a)]$[F(E[p]):F]' \, \leq \, 2$;
\item[(b)]$\mu_p\subset F$.
\end{enumerate}
Then Im$(s_p)\ne 0$.
Moreover, if $F(E[p])=F$,
i.e., if all the $p$-division points are rational over $F$, then
there exist points $P, Q$ of $E(F)/p$ such that
$$
s_p(\langle P, Q\rangle) \ne 0.
$$
In this case, up to replacing $F$ by a finite unramified extension,
we may choose $P$ to be a $p$-power torsion point. \rm

\medskip

\begin{proof}
First consider the case when $E[p]$ is semisimple, so that
$[K:F]'=[K:F] \leq 2$, where $K:=F(E[p])$.

Suppose $K$ is quadratic over $F$. Recall that
over $F$, $C_F$ (resp. $C'_F$) is given by the character
$\chi\nu^{-1}$ (resp. $\nu$), and since $\mu_p\subset F$ and $\nu$
quadratic, we have
$$
H^2(F, C_F^{\otimes 2}) \, \simeq \, H^2(F, \mu_p) \, = \, Br_F[p]
\, \simeq \, \Z/p.
$$
Suppose we have proved the existence of a class $\theta_K$ in
$T_K(E\times E)/p$ such that $s_{p,K}(\theta_K)$ is non-zero, and
this image must be, thanks to Proposition C, in the image of a class
$t_K$ in $Br_K[p]$. Put $\theta:=N_{K/F}(\theta_K) \, \in \, T_F(E
\times E)/p$. Then $s_p(\theta)$ equals
$N_{K/F}(s_{p,K}(\theta_K))$, which is in the image of $t:=
N_{K/F}(t_K)\in Br_F[p]$, which is non-zero because the norm map on
the Brauer group is an isomorphism.

So we may, and we will, assume henceforth in the proof of this
Proposition that $E[p])\subset F$, i.e., $K=F$.

Now let us look at the basic setup carefully. Since $E$ has good
reduction, the N\'eron model is an elliptic curve over $S$.
Moreover, since $E$ has ordinary reduction with $E[p]\subset F$, we
also have
$$
\ScrE[p] \, = \, (\mu_p)_S\oplus\Z/p
$$ as group schemes over $S$
(and as sheaves in $S_{\Flat}$). By the Appendix A.3.2,

\noindent{$(5.1)$}
$$
\begin{CD}\quad E(F)/p\simeq \ScrE(S)/p @>{\partial_F}>>
H_{\fl}^1(S,\ScrE[p])=H_{\fl}^1 (S,\mu_p)\oplus H_{\fl}^1(S,\Z/p),
\end{CD}
$$
where the boundary map $\partial_F$ is an isomorphism. Moreover, the direct sum
on the right of (5.1) is isomorphic to
$$
H_{\fl}^1 (S,\mu_p)\oplus H_{\et}^1(k,\Z/p) \, \simeq \, \ScrO_F^*/p\oplus\Z/p,
$$
thanks to the isomorphism $H_{\fl}^1 (S,\mu_p) \simeq \ScrO_F^*/p$
and the identification of $H_{\fl}^1(S,\Z/p)$ with
$H_{\et}^1(k,\Z/p)\simeq \Z/p$ (\cite{Mi1}, p. 114, Thm.~3.9).
Therefore we have a 1-1 correspondence
$$
\Pbar\longleftrightarrow (\ubar_p,\nbar_p) \leqno{(5.2)}
$$
with $\Pbar\in E(F)/p$, $\ubar_p\in\ScrO_F^*/p$, $\nbar_p\in\Z/p$.
The ordered pair on the right of (5.2) can be viewed as an element
of $H_{\fl}^1 (S,\ScrE[p])$ or of its (isomorphic) image
$\Hbar^1(F,E[p])$ in $H^1(F,E[p])$.

In Galois cohomology, we have the decomposition
$$
H^2(F,E[p]^{\otimes 2}) \, \simeq \, H^2(F,\mu_p^{\otimes 2})\oplus
H^2(F,\mu_p)^{\oplus 2}\oplus H^2(F,\Z/p).\leqno(5.3)
$$
We have a similar one for $H^2_{\fl}(S,\ScrE[p]^{\otimes 2})$.

It is essential to note that by Proposition C,
$$
Im(s_p) \, \hookrightarrow  \, H^2(F,\mu_p^{\otimes 2}) \subset
H^2(F,E[p]^{\otimes 2}),\leqno(5.4-i)
$$
and in addition,
$$
H^2(F,\mu_p^{\otimes 2})\, \simeq \, H^2(F, \mu_p) \, = \, \Br_F[p]
\, \simeq \, \Z/p,\leqno(5.4-ii)
$$
which is implied by the fact that $\mu_p\subset F$ (since it is the
determinant of $E[p])$).

In terms of the decomposition (5.2) from above, we get, for all $P,
Q \in E(F)$,
$$
s_p(\left<P,Q\right>)=\ubar_P\cup\ubar_Q \in Br_F[p] \simeq \Z/p.
\leqno{(5.5)}
$$

\medskip

One knows that (cf. \cite{Lang2}, chapter II, sec. 3, Prop.6)
$$
\vert\ScrO_F^*/p\vert \, = \, \vert\mu_p(F)\vert
p^{[F:\Q_p]}.\leqno(5.6)
$$
Since we have assumed that $F$ contains all the $p$-th roots of
unity, the order of $\ScrO_F^*/p\vert$ is at least $p^2$.

\medskip

\noindent{\bf Claim 5.7} \, \it Let $u, v$ be units in $\ScrO_F$
which are linearly independent in the $\F_p$-vector space $V: =
\ScrO_F^\ast/(\ScrO_F^\ast)^p$. Then $F[u^{1/p}]$ and $F[v^{1/p}]$
are disjoint $p$-extensions of $F$. \rm

\medskip

This is well known, but we give an argument for completeness.
Pick $p$-th roots $\alpha, \beta$ of
$u, v$ respectively. If the extensions are not disjoint, we must have
$\alpha=\sum_{j=0}^{p-1} c_j\beta^j$ in $K=F[v^{1/p}]$, with $\{c_j\}\subset F$.
A generator $\sigma$ of Gal$(K/F$ must send $\beta$ to $w\beta$ for some
$p$-th root of unity $w\ne 1$ (assuming, as we can, that $v$ is not a $p$-th power in $F$),
and moreover, $\sigma$ will send $\alpha$ to $w^i\alpha$ for some $i$.
Thus $\alpha^\sigma$ can be computed in two different ways, resulting in the identity
$\sum_j c_jw^i\beta^j = \sum_j c_jw^j\beta^j$, from which the Claim follows.

\medskip

Consequently, since the dimension of $V$ is at least $2$ and since
$F$ has a unique unramified $p$-extension, we can find
$u\in\ScrO_F^*$ such that $K: = F[u^{1/p}]$ is a ramified
$p$-extension. Fix such a $u$ and let $Q\in E(F)$ be given by
$(\ubar,0)$. By \cite{Serre3}, Prop. 5 (iii) on page 72 (see also
the Remark on page 95), there exists $v\in \ScrO_F^*$ s.t. $v\notin
N_{K/F}(\ScrO_K^*)$. Then by \cite{Tate3}, Prop. 4.3, page 266,
$\{v,u\}\not=0$ in $\Br_F[p]$. Take $P\in E(F)$ such that
$\Pbar\leftrightarrow (\vbar,0)$ then
$s_p(\left<P,Q\right>)=\{v,u\}\in H^2(F,E[p]^{\otimes
2})\subset\Br_F[p]$ and $\{v,u\}\not=0$.

\medskip

It remains to show that we may choose $P$ to be a $p$-power torsion
point after possibly replacing $F$ by a finite unramified extension.
Since $E[p]$ is in $F$, $\mu_p$ is in $F$; recall that $\mu_p(F)
\subset E[p](F)$. So we may pick a non-trivial $p$-th root of unity
$\zeta$ in $F$. Let $m$ be the largest integer such that
$\zeta: = w^{p^{m-1}}$ for some $w\in F^\ast$; in this case $w$ is
in $F^\ast - {F^\ast}^p$. Let $F'/F$ be
the unramified extension of $F$ such that over the corresponding
residual extension, all the $p^m$-torsion points of $\ScrE_s$ are
rational; note however that $\ScrE[p^m]$ need not be in $F$. This
results in the following short exact sequence:
$$
0 \, \to \, \mu_{p^m} \, \to \ScrE[p^m] \, \to \Z/p^m \, \to 0,
$$
leading to the inclusion
$$
\mu_{p^m}(F') \, \subset E[p^m](F').
$$
Since $F'/F$ is unramified, $w$ cannot belong to ${{F'}^\ast}^p$,
and the corresponding point $P$, say, in $E[p^m](F')$ is not in
$pE(F')$. Put
$$
L \, = \, F'(\frac{1}{p}P) \, = \, F'(w^{1/p}),
$$
which is a ramified $p$-extension. So there exists a unit $u$ in
$\ScrO_{F'}$ such that $\{w,u\}$ is not trivial in Br$_{F'}[p]$. Now
let $Q$ be a point in $E(F')$ such that its class in $E(F')/p$ is
given by $(\overline u, 0) \, \in \, \ScrO_{F'}/p \oplus \Z/p$. It
is clear that $\langle P, Q\rangle$ is non-zero in $T_{F'}(A)/p$.

\medskip

Now suppose $E[p]$ is not semisimple. Since we are not in the wild case, we know by Lemma 4.0 that
$K:=F(E[p])$ is an unramified $p$-extension of $F$. By the discussion above, there are points $P, Q \in E(K)/p$ such that $s_{p, K}(\langle P, Q \rangle_K) \ne 0$. Put
$$
\theta: = \, N_{K/F}(\langle P, Q\rangle_K) \, \in T_F(A)/p.
$$

\medskip

\noindent{\bf Claim}: \, $s_p(\theta)$ is non-zero.

\medskip

Since we have
$$
s_p(\theta) \, = \, N_{K/F}\left(s_{p, K}(\langle P, Q \rangle_K)\right),
$$
the Claim is a consequence of the following Lemma, well known to experts.

\medskip

\noindent{\bf Lemma 5.8} \, \it If $K\supset F$ is any finite
extension of $p$-adic fields, then
$N_{K/F}\colon\,\Br_K\to \Br_F$ is an isomorphism. \rm

\medskip

For completeness we give a

\noindent{\it Proof}. \,

The {\it invariant map} $inv_F: \Br(F) \rightarrow \Q/\Z$ is
an isomorphism and moreover (cf. \cite{Serre3}, chap.XIII, Prop.7),
if $n = [K:F]$, the following diagram commutes:
\begin{equation*}
\begin{CD}
{} &\quad {} &\quad \Br_F @>{{\rm Res}_{K/F}}>>
  \Br_K\\
@. @V{{\rm inv}_F}VV @VV{{\rm inv}_K}V\\
&&\Q/\Z @>{n}>> \Q/\Z
\end{CD}
\end{equation*}
As $\Q/\Z$ is divisible, every element in $\Br_K$ is the restriction
of an element in $\Br_F$. Now since
$\Res_{K/F}\colon\,\Br(F)\to\Br(K)$ is already multiplication by $n$
in $\Q/\Z$, the projection formula $N_{K/F} \circ \Res_{K/F} =
[K:F]{\rm Id} = n {\rm Id}$ implies that the norm $N_{K/F}$ on
$\Br_K$ must correspond to the identity on $\Q/\Z$.

\medskip

We are now done with the proof of Proposition E.

\end{proof}

\medskip

\section{Injectivity of $s_p$}
\label{sec6}

\medskip

\noindent{\bf Proposition F} \, \it Let $E$ be an elliptic curve
over a non-archimedean local field $F$ of odd residual
characteristic $p$, such that $E$ has good, ordinary reduction. Then
$s_p$ is injective on $T_F(E\times E)/p$.
\rm

\medskip

In view of Proposition E, we have the following

\medskip

\noindent{\bf Corollary 6.1} \, \it Let $F, E, p$ be as in
Proposition E. Then $T_F(E \times E)/p$ is a cyclic group of
order $p$. Moreover, if $E[p]\subset F$, it even consists of symbols $\langle P,
Q\rangle$, with $P, Q \in E(F)/p$. \rm

\medskip

To prove Proposition F, we will need to consider separately the
cases when $E[p]$ is semisimple and non-semisimple, the latter case being split further
into the wild and unramified subcases.

\bigskip

{\bf Proof of Proposition F in the semisimple case}

\medskip

Again, to prove injectivity, we may replace $F$ by any finite
extension of prime-to-$p$ degree. Since $p$ does not divide
$[F(E[p]):F]$ when $E[p]$ is semisimple, we may assume (in this
case) that all the $p$-torsion points of $E$ are rational over $F$.

\medskip

\begin{Remark}
The injectivity of $s_p$ when $E[p]\subset F$ has been announced
without proof, and in fact for a more general situation, by Raskind
and Spiess \cite{RS}, but the method of their paper is completely
different from ours.
\end{Remark}

\medskip

There are three steps in our proof of injectivity when $E[p]\subset
F$:

\medskip

\noindent {\bf Step I: \, Injectivity of $s_p$ on symbols}

\medskip

Pick any pair of points $P, Q$ in $E(F)$. We have to show that if
$s_p(<P,Q>) =0$, then the symbol $<P,Q>$ lies in $pT_F(A)$.

\medskip

To achieve Step I, it suffices to prove that the condition (a) of
lemma 1.8.1 holds.

\medskip

In the correspondence (5.2), let $\Pbar\leftrightarrow
(\ubar_P,\nbar_P)$ and $\Qbar\leftrightarrow (\ubar_Q,\nbar_Q)$. Put
$K_1=F(\sqrt[p]{u_Q})$ and take $K_2$ to be the unique unramified
extension of $F$ of degree $p$ if $\nbar_Q\not=0$; otherwise take
$K_2=F$. Consider the compositum $K_1K_2$ of $K_1, K_2$, and $K:
=F\left(\frac1p\,Q\right)$, all the fields being viewed as subfields
of $\Fbar$.

\medskip

From (5.1) we get the following commutative diagram:

\medskip

\noindent(6.2)

\begin{equation*}
\begin{matrix}
\ScrE(\ScrO_K)/p &\mapright{\sim} &H_{\fl}^1(\ScrO_K,\ScrE[p])
  &\simeq &\ScrO_K^*/(\ScrO_K^*)^p &\oplus &\Z/p\\
\lmapdown{N_{K/F}}\uparrow{Res} &&\uparrow{Res}&&
\lmapdown{N_{K/F}}\uparrow{Res}
&&\rmapdown{N_{K/F}=\id}\\
\ScrE(\ScrO_F)/p &\mapright{\sim} &H_{\fl}^1
  (\ScrO_F,\ScrE[p]) &\simeq
  &\ScrO_F^*/(\ScrO_F^*)^p &\oplus &\Z/p.
\end{matrix}
\end{equation*}

Here the map {\it Res}$={\it Res}_{K/F}$ on $\ScrE(\ScrO_F)/p$ and $\ScrO_F^\ast/p$
induced by the inclusion $F \hookrightarrow K$ is the obvious
restriction map. However, on $\Z/p$, {\it Res} comes from the
residue fields $\F$ of $F$ and $\F'$ of $K$ via the identifications
$$
H^1_{\rm fl}(S, \Z/p) \simeq H^1(\F, \Z/p) \simeq {\rm Hom}({\rm
Gal}(\overline \F/\F), \Z/p) \simeq \Z/p \leqno(6.3)
$$
and the corresponding one for $K$ and $\F'$. As Gal$(\overline
\F/\F) = \hat \Z$, we see by taking $\F \subset \F' \subset
\overline \F$ with $[\F':\F] = d$ dividing $p$, that Gal$(\overline
\F/\F')$ is the obvious subgroup of Gal$(\overline \F/\F)$
corresponding to $d\hat \Z$; consequently, the map {\it Res} on the
$\Z/p$ summand is multiplication by $d$. Finally note that $N_{K/F}
\circ {\it Res}_{K/F} = [K:F] {\rm id}$.

\medskip

\noindent {\bf Lemma 6.4}\enspace \it With $K, K_1, K_2$ as above,
we have
$$
K \, = \, K_1K_2.
$$
In other words we have the following diagram
\begin{align*}
\vbox{\beginpicture \setcoordinatesystem units <.50cm,.50cm>
\linethickness=1pt \setplotsymbol
({\makebox(0,0)[l]{\tencirc\symbol{'160}}}) \setshadesymbol
({\thinlinefont .}) \setlinear
%
%
\linethickness= 0.500pt \setplotsymbol ({\thinlinefont .}) \plot
11.800 21.969 10.721 20.286 /
%
%
\linethickness= 0.500pt \setplotsymbol ({\thinlinefont .}) \plot
13.968 20.343 12.757 21.933 /
%
%
\linethickness= 0.500pt \setplotsymbol ({\thinlinefont .}) \plot
10.856 19.494 12.067 17.905 /
%
%
\linethickness= 0.500pt \setplotsymbol ({\thinlinefont .}) \plot
14.029 19.592 12.950 17.909 /
%
%
\put{$K=K_1K_2$} [lB] at 11.113 22.257
%
%
\put{\kern-5pt$K_1$} [lB] at  9.970 19.685
%
%
\put{\lower1pt\hbox{$K_2$}} [lB] at 14.192 19.748
%
%
\put{\kern-3pt$F$} [lB] at 12.383 17.367
%
%
\put{$\scriptstyle p$} [lB] at 11.017 18.479
%
%
\put{$\scriptstyle p$} [lB] at 13.621 18.447 \linethickness=0pt
\putrectangle corners at  9.970 22.587 and 14.192 17.285
\endpicture}
&\qquad\raise1.25true cm\hbox{$\begin{matrix} &\text{(If
$\ubar_Q\not=1$,
  $\nbar_Q\not=0$)}\hfill\\
&\text{$K_1\supset F$ ramified}\hfill\\
&\text{$K_2\supset F$ unramified}\hfill
\end{matrix}$}
\end{align*}
\rm

\medskip

\begin{Proof} Put $L = K_1K_2$. Let $\partial_L$ be the boundary map given by the sequence (5.1).

\noindent{a)} \, ${\bf K \subset L}$:  \, Apply the diagram (6.2)
with $L$ instead of $K$. We claim that $\partial_L({\it Res}(Q)) =
0$ and hence $K \subset L$. If $\overline n_Q = 0$, so that $L =
K_1$, we get $\partial_L({\it Res}(Q)) = ({\it Res}(\overline
u_Q),0) = 0$. On the other hand, if $\overline n_Q \ne 0$, then
$K_2$ is the unique unramified $p$-extension of $F$ and the
restriction map on $\Z/p$ is, by the remarks above, given by
multiplication by $p$. So we still have $\partial_L({\it Res}(Q)) =
0$, as claimed.

\noindent{b)} \, ${\bf K \supset L}$:  \, In this case we have
$\partial_K({\it Res}(Q)) = 0$. On the other hand,
$$
\partial_K({\it
Res}(Q)) = ({\it Res}(\overline u_Q),{\it Res}(\overline n_Q)) =
(0,0).
$$
Hence $K \supset K_1$. If $\overline n_Q = 0$, $K_2 = F$ and
we are done. So suppose $\overline n_Q \ne 0$. Then, since ${\it
Res}(\overline n_Q)$ is zero, we see that the residual extension of
$K/F$ must be non-trivial and so must contain the unique unramified
$p$-extension $K_2$ of $F$. Hence $K$ contains $L = K_1K_2$.

\qed
\end{Proof}

\medskip

\noindent{\bf Lemma 6.5} \, \it Let $P, Q$ be in $E(F)/p$ with
coordinates in the sense of (5.2), namely $P = (\overline u_P,
\overline n_P)$ and $Q = (\overline u_Q, \overline n_Q)$. Let $P_0,
Q_0$ be the points in $E(F)/p$ with coordinates $(\overline u_P,
0)$, $(\overline u_Q, 0)$ respectively. Then
$$
\langle P, Q \rangle \, = \, \langle P_0, Q_0 \rangle \, \in
T_F(A)/p.
$$
\rm
\medskip

{\it Proof of Lemma 6.5}. \, Put $P_1 = (0, \overline n_P)$ and $Q_1
= (0, \overline n_Q)$. Then by linearity,
$$
\langle P,Q\rangle \, = \, \langle P_0,Q_0\rangle + \langle P_1,
Q_0\rangle + \langle P_0, Q_1\rangle + \langle P_1, Q_1\rangle.
$$
First note that linearity,
$$
\langle P_1, Q_1\rangle \, = \, \overline n_P\overline n_Q\langle
(0,1), (0,1)\rangle.
$$
It is immediate, since $p \ne 2$, to see that $\langle(0,1),
(0,1)\rangle$ is zero by the skew-symmetry of $\langle ., .
\rangle$. Thus we have, by bi-additivity,
$$
\langle P_1, Q_1\rangle \, = \, 0.
$$

Next we show that in $T_F(A)/p$,
$$
\langle P_0, Q_1\rangle \, = 0 \, = \, \langle P_1 , Q_0\rangle.
$$
We will prove the triviality of $\langle P_0, Q_1\rangle$; the
triviality of $\langle P_1, Q_0\rangle$ will then follow by the
symmetry of the argument. There is nothing to prove if $\overline
n_Q = 0$, so we may (and we will) assume that $\overline n_Q \ne 0$.
Then it corresponds to the unramified $p$-extension $M =
F\left(\frac{1}{p}Q_1\right)$ of $F$. It is known that every unit in
$F^\ast$ is the norm of a unit in $M^\ast/p$. This proves that the
point $P_0$ in $E(F)/p$ is a norm from $M$. Thus $\langle P_0,
Q_1\rangle$ is zero by Lemma 1.8.1. Putting everything together, we
get
$$
\langle P, Q\rangle \, = \, \langle P_0, Q_0 \rangle
$$
as asserted in the Lemma.

\qed

\medskip

\noindent{\it Proof of Step I}:

\medskip

Let $P, Q \in E(F)/p$ be such that $s_p(\left<P,Q\right>)=0$. We
have to show that $\langle P, Q \rangle = 0$ in $T_F(A)/p$. By Lemma
6.5 we may assume that $\overline n_P=\overline n_Q=0$. So it follows that the
element $\{\ubar_P,\ubar_Q\}:=\ubar_P \cup \ubar_Q$ is zero in the
Brauer group $\Br_F[p]$. Then putting $K_1 =
K\left(u_Q^{1/p}\right)$, we have by \cite{Tate3}, Prop. 4.3, we
have $u_P=N_{K_1/F}(u'_1)$ with $u'_1\in K_1^*$ and where
$u_P\in\ScrO_F^*$ with image $\ubar_P$. It follows that $u'_1$ must
also be a unit in $\ScrO_{K_1}$. Now with the notations introduced
at the beginning of {\it Step I} we have, since $\overline n_Q =0$,
that $K=K_1$ and so $[K:F] = p$. In the diagram (6.2), take in the
upper right corner the element $(\overline{u'},0)$. Then we get an
element $\overline{P'}\in E(K)/p$ such that
$N_{K/F}\overline{P'}\equiv P$ $(\mod\,p)$. Hence condition (a) of
Lemma 1.8.1 holds, yielding Step I.

\qed

\bigskip

\noindent {\bf Step II}

\medskip

Put
$$
{\ST}_{F,p}(A) \, = \, {\rm Im}\left(\ST_F(A) \, \rightarrow \,
T_F(A)/p\right)
$$

\medskip

\noindent{\bf Lemma 6.6} \, \it  $s_p$ is injective on
${\ST}_{F,p}(A)$.
\rm

\medskip

\begin{Proof}
This follows from

\medskip

\smallskip\noindent
{\bf Claim:}\enspace If $V$ is an $\F_p$-vector space with an
alternating bilinear form $[\quad,\quad]\colon\, V\times V\to W$,
with $W$ a $1$-dimensional $\F_p$-vector space, then the
conditions (i) and (ii) of \cite{Tate3}, p. 266, are satisfied,
which in our present setting read as follows :
\begin{enumerate}
\item[(i)] Given $a,b, c, d$ in $V$ such that $[a,b] = [c,d]$, then
there exist elements $x$ and $y$ in $V$ such that
$$
[a,b] = [x,b] =[x,y] = [c,y] = [c,d].
$$
\item[(ii)] Given $a_1,a_2, b_1, b_2$  in $V$, there exist $c_1 ,c_2$ and
$d$ in $V$ such that
$$
[a_1,b_1] = [c_1,d ], \, \, {\rm and} \, \,  [a_2,b_2] = [c_2,d] .
$$
\end{enumerate}

\medskip

\begin{Proof} Recall that we are assuming in this paper that $p$ is odd.
The Claim is just Proposition 4.5 of \cite{Tate3}, p. 267. The proof goes
through verbatim. In our case we apply this to $V=E(F)/p,
W=\Br_F[p]$ and $[\Pbar,\Qbar]=s_p(\left<P,Q\right>) =\ubar_P\cup
\ubar_Q$. By (\cite{Tate3}, Corollary on p. 266), the $s_p$ is injective
on ${\ST}_{F,p}(A)$.
\end{Proof}

\bigskip

\noindent
{\bf Step III: \, Injectivity of $s_p$ in the general case} (but
still assuming $E[p]\subset F$).

\medskip

For this we shall appeal to Lemma 5.8.

\medskip

\noindent
{\it Proof of Step III}

\medskip

By Step II, we see from Lemma 1.7.1 that it suffices to prove the
following result, which may be of independent interest.

\medskip

\noindent{\bf Proposition 6.7} \, \it Let $K \supset F$ be a finite
extension, and $E[p] \subset F$. Then $N_{K/F}({\ST}_{K,p}(A))$ is a
subset of ${\ST}_{F,p}(A)$ and hence
$$
{\ST}_{F,p}(A) \, = \, T_F(A)/p.
$$
In other words, $T_F(A)$ is generated by symbols modulo $pT_F(A)$.
\rm

\medskip

It is an open question as to whether $T_F(A)$ is different from
${\ST}_F(A)$, though many expect it to be so.

\medskip

{\it Proof of Proposition 6.7}. \, We have to prove that if $P',Q'
\in E(K)$, then the norm to $F$ of $\langle P', Q' \rangle$ is
mapped into ${\ST}_{F,p}(A)$. Since we have
$$
{\rm Im} s_p({\ST}_F(A)) \, = \, {\rm Im} s_p(T_F(A)/p \, \simeq \,
\Z/p,
$$
the assertion is a consequence of the following

\medskip

\noindent{\bf Lemma 6.8} \, \it  If $[K\colon\,F]=n$ and $P',Q'\in
E(K)$, $P,Q\in E(F)$ are such that
$$
s_{p,F}(N_{K/F}\left<P',Q'\right>) \, = \,
s_{p,F}(\left<P,Q\right>),
$$
then, assuming that all the $p$-torsion in $E(\overline F)$ is
$F$-rational,
$$
N_{K/F}(\left<P',Q'\right>)\equiv \left<P,Q\right> \, (\mod\,
p\,T_F(A))
$$.
\rm

\medskip

{\it Proof of Lemma 6.8}. \, We start with a few simple sublemmas.

\medskip

\noindent{\bf Sublemma 6.9} \, \it If $K/F$ is unramified, then
Lemma 6.8 is true. \rm

\medskip

{\it Proof}. \, Indeed, in this case, every point in $E(F)$ is the
norm of a point in $E(K)$ (cf. \cite{Mazur}, Corollary 4.4). So we
can write $P \equiv N_{K/F}(P_1) (\bmod pE(K))$ with $P_1 \in E(K)$,
from which it follows that
$$
N_{K/F}(\langle P_1, Q\rangle) \, \equiv \, \langle P,Q\rangle \,
(\bmod pT_F(A)).
$$
On the other hand, by the remark above (in the beginning of Step
III) about the norm map on the Brauer group, we have
$$
s_{p,K}(\langle P_1,Q\rangle) \, = \, s_{p,K}(\langle P',Q'\rangle)
$$
and applying Step II to $\langle P_1, Q\rangle$ and $\langle P',
Q'\rangle$ for $K$, we are done.

\qed

\medskip

\noindent{\bf Subemma 6.10} \, \it If $K/F = n$ with $p \nmid n$,
then Lemma 6.8 is true. \rm

\medskip

{\it Proof}. \, Indeed, as the cokernel of $N_{K/F}$ is annihilated
by $n$ which is prime to $p$, the norm map on $E(K)/p$ is surjective
onto $E(F)/p$. Let $P_1 \in E(K)$ satisfy $P \equiv N_{K/F}(P_1)
(\bmod pE(F))$. Then by the projection formula,
$$
\langle P, Q\rangle \, \equiv \, N_{K/F}(\langle P_1, Q\rangle) \,
(\bmod pT_F(A)). \leqno(\ast)
$$
Since $s_{p,F} \circ N_{K/F} = N_{K/F} \circ s_{p,K}$ and since
$N_{K/F}$ is non-trivial, and hence an isomorphism, on the
one-dimensional $\F_p$-space $\Br_K[p]$, we see that
$$
s_{p,K}(\langle P', Q'\rangle) \, = \, s_{p,K}(\langle P_1,
Q\rangle).
$$
Applying Lemma 6.6 for $K$ we see that $\langle P', Q'\rangle)$
equals $\langle P_1, Q\rangle$ in $T_K(A)/p$. Now we are done by
$(\ast)$.

\qed

\medskip

\noindent{\bf Sublemma 6.11} \, \it If $K\supset F_1 \supset F$ and
if 6.8 is true for both the pairs $K/F_1$ and $F_1/F$, then it is
true for $K/F$.\rm

\medskip

{\it Proof}. \, Let $P,Q \in E(F)$ and $P', Q' \in E(K)$ be such
that
$$
s_{p,F}(N_{K/F}\left<P',Q'\right>) \, = \,
s_{p,F}(\left<P,Q\right>). \leqno(a)
$$
Applying Proposition E with $F_1$ in the place of $F$, we get points
$P_1, Q_1$ in $E(F_1)$ such that
$$
s_{p,F_1}(\langle P_1,Q_1\rangle) \, = \, N_{K/F_1}(s_{p,K}(\langle
P',Q'\rangle)). \leqno(b)
$$
Since the right hand side is the same as
$s_{p,F_1}(N_{K/F_1}(\langle P',Q'\rangle))$, we may apply 6.8 to
the pair $K/F_1$ to conclude that
$$
\langle P_1,Q_1\rangle \, = \, N_{K/F_1}(\langle P',Q'\rangle).
\leqno(c)
$$
Applying $N_{F_1/F}$ to both sides of $(b)$, using the facts that
the norm map commutes with$s_p$ and that $N_{K/F} = N_{K/F_1} \circ
N_{F_1/F}$, and appealing to $(a)$, we get
$$
s_{p,F}(N_{F_1/F}\left<P_1,Q_1\right>) \, = \, s_{p,F}(\langle
P,Q\rangle). \leqno(d)
$$
Applying (6.8) to $F_1/F$, we then get
$$
\langle P,Q\rangle \, = \, N_{F_1/F}(\langle P_1,Q_1\rangle).
\leqno(e)
$$
The assertion of the Sublemma now follows by combining $(c)$ and
$(e)$.

\qed

\medskip

\noindent{\bf Sublemma 6.12} \, \it It suffices to prove 6.8 for all
finite Galois extensions $K'/F'$ with $E[p] \subset F'$ and
$[F':\Q_p] < \infty$. \rm

\medskip

{\it Proof}. Assume 6.8 for all finite Galois extensions $K'/F'$ as
above. Let $K/F$ be a finite, non-normal extension with Galois
closure $L$, and let $E[p] \subset F$. Suppose $P, Q \in E(F)$ and
$P', Q' \in E(K$ satisfy the hypothesis of 6.8. We may use the
surjectivity of $N_{L/K}: \Br_L \to \Br_K$ and Proposition E (over
$L$) to deduce the existence of points $P'', Q'' \in E(L)$ such that
$$
s_{p,K}(N_{L/K}\left<P'',Q''\right>) \, = \,
s_{p,K}(\left<P',Q'\right>).
$$
As $L/K$ is Galois, we have by hypothesis,
$$
N_{L/K}(\left<P'',Q''\right>) \, = \, \left<P',Q'\right>. \leqno(i)
$$
By construction, we also have
$$
s_{p,F}(N_{L/F}\left<P'',Q''\right>) \, = \,
s_{p,F}(\left<P,Q\right>).
$$
As $L/F$ is Galois, we have
$$
N_{L/F}(\left<P'',Q''\right>) \, = \, \left<P,Q\right>. \leqno(ii)
$$
The assertion now follows by applying $N_{K/F}$ to both sides of (i)
and comparing with (ii).

\qed

\medskip

Thanks to this last sublemma, we may assume that $K/F$ is Galois.
Appealing to the previous three sublemmas, we may assume that we are
in the following {\bf key case}:

\medskip

\noindent${\bf ({\mathcal K})}$ \, \it $K/F$ is a totally ramified,
cyclic extension of degree $p$, with $E[p] \subset F$ \rm

\medskip

So it suffices to prove the following

\medskip

\noindent{\bf Lemma 6.13} \, \it 6.8 holds in the key case
$({\mathcal K})$. \rm

\medskip

{\it Proof}. \, Suppose $K/F$ is a cyclic, ramified $p$-extension
with $E[p] \subset F$, and let $P, Q \in E(F)$, $P', Q' \in E(K)$
satisfy the hypothesis of 6.8. Since $\mu_p = \det(E[p])$, the
$p$-th roots of unity are  in $F$ and so
$$
W: = F^\ast/p \, = \, H^1(F, \mu_p) \, \simeq \, H^1(F, \Z/p)
\, = \, {\rm Hom}({\rm Gal}(\overline F/F), \Z/p).
$$
In other words, lines in $W$ correspond to cyclic $p$-extensions of
$F$, and we can write $K = F\left(u^{1/p}\right)$ for some $u \in
\ScrO_F^\ast$. For every $w \in \ScrO_F^\ast$, let $\overline w$
denote its image in $W$. Put $m = \dim_{\F_p} W$. Then $m \geq 3$
by \cite{Lang2}, chapter II, sec.~3, Proposition 6.
Let $W_1$ denote the line spanned in $W$ by the unique unramified
$p$-extension of $F$. Since $m > 2$, we can find some $v \in
\ScrO_F^\ast$ such that $\overline v$ is not in the linear span of
$W_1$ and $\overline u$. Using the {\it Claim} stated in the proof
of Proposition E, we see that $L: = F\left(v^{1/p}\right)$ is
linearly disjoint from $K$ over $F$. Both $K$ and $L$ are totally
ramified $p$-extensions of $F$, and so is the $(p, p)$-extension
$KL$. Then $KL/K$ is a cyclic, ramified $p$-extension, and by local
class field theory (\cite{Serre3}, chap. V, $\sec 3$), there exists
$y \in \ScrO_K^\ast$ not lying in $N_{KL/K}((KL)^\ast)$. Hence
$\overline v \cup \overline y$ is non-zero in $\Br_K[p]$ (see
\cite{Tate3}, Prop.4.3, p.266, for example). Let $P_1 \in E(F)$ and
$Q_1 \in E(K)$ be such that in the sense of (5.2) we have
$$
\overline P_1 \leftrightarrow (\overline v, 0), \, \, \overline Q_1
\leftrightarrow (\overline y, 0).
$$
Then
$$
s_{p,K}(\langle P_1,Q_1\rangle) \, = \, \overline v \cup \overline y
\, \ne \, 0.
$$
Since $\Br_K[p]$ is one-dimensional (over $\F_p$), we may replace
$Q_1$ by a multiple and assume that
$$
s_{p,K}(\langle P_1,Q_1\rangle) \, = \, s_{p,K}(\langle
P',Q'\rangle).
$$
Hence
$$
\langle P_1,Q_1\rangle \, \equiv \, \langle P',Q'\rangle \, \bmod
pT_K(A)
$$
by {\it Step II} for $K$. So we get by the hypothesis,
$$
s_{p,F}(N_{K/F}\left<P_1,Q_1\right>) \, = \,
s_{p,F}(\left<P,Q\right>).
$$
But $P_1$ is by construction in $E(F)$, and so the projection
formula says that
$$
N_{K/F}\left<P_1,Q_1\right> \, = \, \langle P_1,
N_{K/F}(Q_1)\rangle.
$$
Now we are done by applying {\it Step II} to $F$.

\end{Proof}

\qed

\medskip

Now we have completed the proof of Proposition F when $E[p]$
is semisimple.

\vskip 0.2in

{\bf Proof of Proposition F in the non-semisimple, unramified case}

\medskip

In this case $K:=F(E[p])$ is an unramified $p$-extension of $F$. By Mazur
(\cite{Mazur}), we know that as $E/F$ is ordinary, the norm map from $E(K)$ to $E(F)$ is surjective. This  shows the injectivity on symbols, giving Step I in this case. Indeed, if $s_p(\langle P, Q\rangle)=0$ for some $P, Q\in E(F)/p$, we may pick an $R\in E(K)$ with norm $Q$ and obtain
$$
N_{K/F}(s_p(\langle P, R\rangle_K)) \, = \, s_p(\langle P, N_{K/F}(R)\rangle) \, = \, 0.
$$
On the other hand, the norm map on the Brauer group is injective (cf. Lemma 5.8). This forces $s_p(\langle P, R\rangle_K)$ to be non-zero. As $E[p]$ is trivial over $K$, the injectivity in the semisimple situation gives the triviality of $\langle P, R\rangle_K = 0$.

The proof of injectivity on $ST_{F,p}(A)$ (Step II) is the same as in the semisimple case.

The proof of injectivity on $T_F(A)/p$ (Step III) is an immediate consequence of the following Lemma (and the proof in the semisimple case).

\medskip

\noindent{\bf Lemma 6.14} \, \it Let $\theta\in T_F(A)/p$. Then $\exists \, \tilde\theta\in T_{K}(A)/p$ such that $\theta = N_{K/F}(\tilde\theta)$. \rm

\medskip

{\it Proof of Lemma}. \, In view of Lemma 1.7.1, it suffices to show, for any finite extension $L/F$ and points $P, Q \in E(L)/p$, that $\beta:= N_{L/F}(\langle P, Q\rangle_L)$ is a norm from $K$. This is obvious if $L$ contains $K$. So let $L$ not contain $K$, and hence is linearly disjoint from $K$ over $F$ (as $K/F$ has prime degree). Consider the compositum $M=LK$, which is an unramified $p$-extension of $L$. By Mazur (\cite{Mazur}), we can write $Q= N_{M/L}(R)$, for some $R \in E(M)/p$. Put
$$
\tilde\beta: = \, N_{M/K}(\langle P, R\rangle_M) \, \in T_K(A)/p.
$$
Then, using $N_{L/F}\circ N_{M/L} = N_{M/F} = N_{K/F}\circ N_{M/K}$, we obtain the desired conclusion
$$
\beta = N_{L/F}\left(\langle P, N_{M/L}(R)\rangle_L\right) = N_{M/F}\left(\langle P, R\rangle_M\right) =  N_{K/F}(\tilde\theta).
$$

\qed

This finishes the proof when $E[p]$ is non-semisimple and unramified.

\vskip 0.2in

{\bf Proof of Proposition F in the non-semisimple, wild case}

\medskip

Here $K/F$ is a wildly ramified extension, so we cannot reduce to the case
$E[p]\subset F$. Still, we may, as in section 4, assume that
$\mu_p\subset F$ and that $\nu=1$, which implies that the
$p$-torsion points of the special fibre $\ScrE_s$ are rational over
the residue field $\F_q$ of $F$. We have in effect (4.2) through
(4.4).

Consider the commutative diagram (4.4). We can write
$$
\ScrO_F^\ast/p \, = \, H^1_{\fl}(S, \mu_p) \, \simeq \, \overline
H^1(F, \mu_p) \, = \, X + Y_S,
$$
with $X \simeq \Z/p$ and $Y_S$ some complement. Furthermore, we have
(cf. \cite{Mi1}, Theorem 3.9, p.114)
$$
H^1_{\fl}(S, \Z/p) \, \simeq \, H^1_{\rm et}(k, \Z/p) \, = \, \Z/p.
$$

Hence
$$
H^1_{\fl}(S, \ScrE[p]) \, = \, Z_S \oplus \Z/p, \quad {\rm where}
Z_S = \psi_S(Y_S). \leqno(6.15)
$$
(Note that the surjectivity of the right map on the top row of the
diagram (4.4) comes from the fact that $H^2_{\fl}(S, \mu_p)$ is $0$
(\cite{Mi2}, chapter III, Lemma 1.1). Recall (cf. (6.2)) that
$E(F)/p$ is isomorphic to $H^1_{\fl}(S, \ScrE[p])$. So by (6.15), we
have a bijective correspondence
$$
P \, \longleftrightarrow \, (\tilde u_P, \overline n_P) \leqno(6.16)
$$
where $P$ runs over points in $E(F)/p$. Compare this with (5.2).

\medskip

\noindent{\bf Lemma 6.17} \, \it Fix any odd prime $p$. Let $K$ be an
arbitrary finite extension of $F$ where $E[p]$ remains
non-semisimple as a ${G}_K$-module. Suppose $\langle(u',0),
(u'',0)\rangle \, = \, 0$ in $T_K(A)/p$ for all $u', u'' \in
\ScrO_K^\ast/p$. Then $\langle P, Q\rangle \, = \, 0$ for all $P, Q
\in E(F)/p$. \rm

\medskip

{\it Proof of Lemma}. \, Let $P = (u_P, n_P), Q = (u_Q, n_Q)$ be in
$E(F)/p$. Then by the bilinearity and skew-symmetry of $\langle . ,
. \rangle$, the fact that $\langle 1, 1 \rangle = 0$, and also the
hypothesis of the Lemma (applied to $K = F$), it suffices to show
that
$$
\langle (u_P, 0) , (0, n_Q) \rangle \, = \, \langle (u_Q, 0) , (0,
n_P) \rangle \, = \, 0.
$$
It suffices to prove the triviality of the first one, as the
argument is identical for the second. Let $L$ denote the unique
unramified $p$-extension of $F$. Then there exists $u' \in
\ScrO_L^\ast/p$ such that $u_P = N_{L/F}(u')$. Since $\langle
(u_P,0), (0,n_Q) \rangle$ equals $N_{L/F}(\langle (u_P, 0), {\rm
Res}_{L/F}((0,n_Q))\rangle)$ by the projection formula, it suffices
to show that
$$
\langle (u', 0), {\rm Res}_{L/F}((0,n_Q))\rangle \, = \, 0.
$$
Now recall (cf. the discussion around (6.3)) that the restriction
map $H^1(\ScrO_F, \Z/p) \rightarrow H^1(\ScrO_L, \Z/p)$ is
zero. This implies that ${\rm Res}_{L/F}((0,n_Q)) = (u'',0)$ for
some $u'' \in \ScrO_L^\ast/p$. (We cannot claim that this
restriction is $(0,0)$ in $E(L)/p$ because the splitting of the
surjection $H^1(\ScrO_K, \ScrE[p]) \to H^1(\ScrO_K, \Z/p)$ is not
canonical when $E[p]$ is non-semisimple over $K$. This point is what
makes this Lemma delicate.) Thus we have only to check that
$\langle(u',0), (u'',0)\rangle \, = \, 0$ in $T_L(A)/p$. This is a
consequence of the hypothesis of the Lemma, which we can apply to $K
= L$ because $E[p]$ remains non-semisimple over any unramified
extension of $F$. Done.

\medskip

As in the semisimple case, there are three steps in the proof.

\bigskip

\noindent{\bf Step I} \, {\bf Injectivity on symbols}:

\medskip

We shall use the following terminology. For $u \in \ScrO_F^\ast$,
write $\overline u$ and $\tilde u$ for its respective images in
$\ScrO_F^\ast/p$ and $Y$, seen as a quotient of $H^1_{\fl}(S,
\mu_p)/X$. We will also denote by $\tilde u$ the corresponding
element in $Z = \psi_S(Y_S) \subset H^1_{\fl}(S, \ScrE[p])$, which
should not cause any confusion as $Y$ is isomorphic to $Z$. We use
similar notation in $\overline H^1(F, -)$. For $\overline u$ in
$\ScrO_F^\ast/p$, we denote by $P_{\overline u}$ the element in
$E(F)/p$ corresponding to the pair $(\tilde u,0)$ given by (6.16).

As in the proof of Proposition D under the diagram (4.4), pick a non-zero element $e \in
X \subset \ScrO_F^\ast/p$. By definition, $\tilde e = 0$, so that
$P_e = (\tilde e,0)$ is the zero element of $E(F)/p$. As we have
seen in the proof of Proposition D, there exists a $v$ in
$\ScrO_F^\ast/p$ such that $[\overline e, \overline v]: = \overline
e \cup \overline v$ is non-zero in Br$_F[p]$.

Now let $\overline x, \overline y \in \ScrO_F^\ast/p$ be such that
$s_p(\langle P_{\overline x}, P_{\overline y}\rangle) \, = \, 0$. We
have to show that
$$
\langle P_{\overline x}, P_{\overline y}\rangle \, = \, 0 \, \in
T_F(A)/p. \leqno(6.18)
$$
There are two cases to consider.

\medskip

\noindent{\bf Case (a)} \, $[\overline x, \overline y] = 0 \in
\Br_F[p]$:

Then $\overline x$ is a norm from $K = F(y^{1/p})$. This implies, as
in the proof in the semisimple case, that $P_{\overline x}$ is a
norm from $E(K)/p$, and by Lemma 1.8.1, we then have (6.18).

\medskip

\noindent{\bf Case (b)} \, $[\overline x, \overline y] \ne 0$:

Since $\Br_F[p] \simeq \Z/p$, we may, after modifying by a scalar,
assume that $[\overline x, \overline y] = [\overline e, \overline
v]$ in $\Br_F[p]$. Then by Tate (\cite{Tate3}, page 266, conditions
(i), (ii)), there are elements $\overline a, \overline b \in
\ScrO_F^\ast/p$ such that
$$
[\overline x, \overline y] = [\overline x, \overline b] = [\overline
a, \overline b] = [\overline a, \overline v] = [\overline e,
\overline v].
$$
\noindent{\bf Claim 6.19} \, \it For each pair of neighbors in this
sequence, the corresponding symbols are equal in $T_F(A)/p$. \rm

\medskip

{\it Proof of Claim 6.19}. \, From $[\overline x, \overline y] =
[\overline x, \overline b]$ we have, by linearity, that $[\overline
x, \overline y\overline b^{-1}] = 0$, hence $x$ is a norm from $L: =
F((yb^{-1})^{1/p})$. Hence $P_{\overline x}$ is a norm from from
$E(L)/p$, and so we get $\langle P_{\overline x},
P_{\overline{yb^{-1}}}\rangle = 0$ in $T_F(A)/p$. Consequently, now
by the linearity of $\langle \cdot , \cdot\rangle$, $\langle
P_{\overline x}, P_{\overline y}\rangle \, = \,  \langle
P_{\overline x}, P_{\overline b}\rangle$. The remaining assertions
of the Claim are proved in exactly the same way.

\qed

\medskip

This Claim finishes the proof of Step I since $\langle P_{\overline
e}, P_{\overline v}\rangle \, = \, 0$, which holds because
$P_{\overline e} = 0$ in $E(F)/p$.
\qed

\medskip

\noindent{\bf Step II} \, {\bf Injectivity on ${\bf ST_{F,p}(A)
\subset T_F(A)/p}$} :

\medskip

{\it Proof}. \, We have just seen in  Step I that each of the symbol
in $T_F(A)/p$ is zero. Since $ST_{F,p}(A)$ is generated by such
symbols, it is identically zero. Done in this case.

\medskip

\noindent{\bf Step III} \, {\bf Injectivity on ${\bf T_F(A)/p}$} :

\medskip

{\it Proof}. \, Since $T_F(A)$ is generated by norms of symbols from
finite extensions of $F$, it suffices to show the following for {\it
every finite extension $L/F$ and points $P, Q \in E(L)$}:
$$
N_{L/F}(\langle P, Q\rangle \, = \, 0 \, \in T_F(A)/p. \leqno(6.20)
$$
Fix an arbitrary finite extension $L/F$ and put
$$
F_1 \, = \, F(E[p]).
$$
There are again two cases.

\medskip

\noindent{\bf Case (a)} \, $L \supset F_1$:

\medskip

Here we are in the situation where $E[p] \subset L$. In particular,
$s_{p,L}$ is injective.

In the correspondence of (5.2), let $P, Q \in E(L)/p$ correspond to
$\overline u_P, \overline u_Q \in \ScrO_L^\ast/p$ respectively. Put
$$
t : \, = \, s_{p,L}(\langle P, Q\rangle) \, = \, [\overline u_P,
\overline u_Q] \, \in \Br_F[p].
$$
If $t = 0$, we have $\langle P, Q\rangle = 0$ in $T_L(A)/p$ (as
$E[p]\subset L$), and we are done.

So let $t \ne 0$. Now let $e$, as before, be the $\F_p$-generator of
$X \subset \ScrO_F^\ast/p$, where $X$ is the image of $\Z/p$
encountered in the proof of Proposition D. The diagram (4.4)
shows that $F_1 = F(e^{1/p})$. In particular, $F_1/F$ is a ramified
$p$-extension. Since $p \ne 2$, $\ScrO_F^\ast/p$ has dimension at
least $3$, and so we can take $v \in \ScrO_F^\ast/p$ such that the
space spanned by $\overline e$ and $\overline v$ in $\ScrO_F^\ast/p$
has dimension $2$ and does not contain the line corresponding to the
unique unramified $p$-extension $M$of $F$. Put
$$
F_2 \, = \, F(v^{1/p}) \quad {\rm and} \quad K \, = \, F_1F_2 \,
\subset \, \overline F.
$$
Then $K/F$ is totally ramified with $[K:F] = p^2$. Hence $K/F_1$ is
still a ramified $p$-extension. By the local class field theory
(\cite{Serre3}, chapter V, section 3), we can choose $y \in
\ScrO_{F_1}^\ast/p$ such that $y \notin N_{K/F_1}(K^\ast)$, so that
$[\overline v, \overline y] \ne 0$ in Br$_F[p]$. By linearity, we
can replace $y$ by a power such that $[\overline v, \overline y] = t
= [\overline u_P, \overline u_Q]$. Now  we have
$$
\langle P_{\overline v}, P_{\overline y}\rangle_{L} \, = \, \langle
P, Q\rangle_{L} \, \in \, T_L(A)/p, \leqno(6.21)
$$
where $\langle \cdot, \cdot \rangle_{L}$ denotes the pairing over
$L$. (This makes sense here because of the hypothesis $F_1 \subset
L$.) Applying the projection formula to $L/F_1$ and $F_1/F$, using
the facts that $P_{\overline v} \in E(F)/p$ and $P_{\overline y} \in
E(F_1)/p$,
$$
N_{L/F_1}(\langle P_{\overline v}, P_{\overline y}\rangle_{L}) \, =
\, [L:F_1]\langle P_{\overline v}, P_{\overline y}\rangle_{F_1}
\leqno(6.22)
$$
and
$$
N_{F_1/F}(\langle P_{\overline v}, P_{\overline y}\rangle_{F_1}) \,
= \, \langle P_{\overline v}, N_{F_1/F}(P_{\overline y})\rangle.
$$
Now since $N_{L/F} = N_{L/F_1} \circ N_{F_1/F}$ and since we have
shown that the symbols in $T_F(A)/p$ are all trivial, we get the
triviality of $N_{L/F}(\langle P, Q \rangle)$ by combining (6.21) and
(6.22).

\medskip

\noindent{\bf Case (b)} \, $L  \not\supset F_1$:

\medskip

In this case $E[p]$ is not semisimple over $L$. Hence by Step II, we
have $\langle P, Q\rangle_L = 0$ in $T_L(A)/p$. So (6.20) holds.
Done.

This finishes the proof of Step III, and hence of Proposition F.

\qed

\medskip

\section{Remarks on the case of multiplicative reduction}
\label{sec7}

\medskip

Here one has the following version of Theorem A, which we will need
in the sequel giving certain global applications:

\medskip

\noindent{\bf Theorem G} \quad \it Let $F$ be a non-archimedean
local field of characteristic zero with residue field $\F_{q}$,
$q=p^r$. Suppose $E/F$ is an elliptic curve over $F$, which has
multiplicative reduction. Then for any prime $\ell\ne 2$, possibly
with $\ell =p$, the following hold:
\begin{enumerate}
\item[{(a)}] \, $s_\ell$ is injective, with image of $\F_\ell$-dimension $\leq
1$.
\item[{(b)}] \, Im$(s_\ell)=0$ if $\mu_\ell\not\subset F$.
\end{enumerate}
\rm

\medskip

This theorem may be proved by the methods of the previous sections,
but in the interest of brevity, we will just show how it can be
derived, in effect, from the results of T.~Yamazaki (\cite{Y}).

\medskip

The Galois representation $\rho_F$ on $E[\ell]$ is reducible as in
the ordinary case, giving the short exact sequence (cf.
\cite{Serre2}, Appendix)
$$
0\to C_F \to E[\ell]\to C'_F \to 0,\leqno(7.1)
$$
where $C'_F$ is given by an unramified character $\nu$ of order
dividing $2$; $E$ is split multiplicative iff $\nu=1$. And $C_F$ is
given by $\chi\nu$ ($=\chi\nu^{-1}$ as $\nu^2=1$), where $\chi$ is the
mod $\ell$ cyclotomic character given by the Galois action on $\mu_\ell$.

Using (7.1), we get a homomorphism (as in the ordinary case)
$$
\gamma_F: \, H^2(F, C_F^{\otimes 2}) \, \rightarrow \, H^2(F,
E[\ell]^{\otimes 2}).\leqno(7.2)
$$
Here is the analogue of Proposition C:

\medskip

\noindent{\bf Proposition 7.3} \, \it We have
$$
Im(s_\ell) \, \subset \, Im(\gamma_F).
$$
Furthermore, the image of $s_\ell$ is at most one-dimensional as
asserted in Theorem G. \rm

\medskip

{\it Proof of Proposition 7.3} \, As $\nu^2=1$ and $\ell\ne 2$, we
may replace $F$ by the (at most quadratic) extension $F(\nu)$ and
assume that $\nu=1$, which implies that $C_F=\mu_\ell$ and
$C'_F=\Z/\ell$. We get the injective maps
$$
\psi: H^1(F,\mu_\ell)\to H^1(F, E[\ell])
$$
and
$$
\delta: E(F)/\ell \to H^1(F, E[\ell]).
$$

\medskip

\noindent{\bf Claim 7.4} \, \it $Im(\delta) \subset Im(\psi)$. \rm

\medskip

Indeed, since $E(F)$ is the Tate curve $F^\ast/q^\Z$ for some $q\in
F^\ast$ (cf. \cite{Serre2}, Appendix), we have the uniformizing map
$$
\varphi: F^\ast \to E(F).
$$
Given $P \in E(F)$, pick an $x\in F^\ast$ such that $P=\varphi(x)$.
Now $\delta(P)$ is represented by the $1$-cocycle $\sigma\to
\sigma(\frac1\ell P)-\frac1\ell P$, for all $\sigma\in G_F$. Let
$y\in {\overline F}^\ast$ be such that $y^\ell=x$. Put $K=F(y)$. One
knows that there is a natural uniformization map (over $K$)
$\varphi_K:K^\ast\to E(K)$, which agrees with $\varphi$ on $F^\ast$.
Then we may take $\varphi_K(y)$ to be $\frac1\ell P$. It follows
that the $1$-cocycle above sends $\sigma$ to
$\frac{y^\sigma}{y}=\zeta_\sigma$, for an $\ell$-th root of unity
$\zeta_\sigma$. In other words, the class of this cocycle is in the
image of $\psi$ in the long exact sequence in cohomology:
$$
0\to\mu_\ell(F)\to E[\ell](F)\to\Z/\ell\to
H^1(F,\mu_\ell)\mapright{\psi} H^1(F,E[\ell])\to \Z/\ell\to \dots,
$$
namely $\psi((\zeta_\sigma)) = \delta(P)(\sigma)$, with
$(\zeta_\sigma)\in H^1(F, \mu_\ell)$. Hence the Claim.

\medskip

Thanks to the Claim, the first part of the Proposition follows
because Im$(s_\ell)$ is obtained via the cup product of classes in
Im$(\delta)$, and Im$(\gamma_F)$ is obtained via the cup product of
classes in Im$(\psi)$.

The second part also follows because $H^2(F,\mu_\ell^{\otimes 2})$
is isomorphic to (the dual of) $H^0(F, \Z/\ell(-1))$, which is at most
one-dimensional.
\qed

\medskip

The injectivity of $s_\ell$ has been proved in Theorem 4.3 of
\cite{Y} in the split multiplicative case, and this saves a long
argument we can give analogous to our proof in the ordinary
situation. Injectivity also follows for the general case because we
can, since $\ell\ne 2$, make a quadratic base change to $F(\nu)$
when $\nu\ne 1$. When combined with Proposition 7.3, we get part (a)
of Theorem G.

\medskip

To prove part (b) of Theorem G, just observe that, thanks to
Proposition 7.3, it suffices to prove that when $\mu_\ell\notin F$,
we have $H^2(F, \mu_\ell^{\otimes 2})=0$. This is clear because the
one-dimensional Galois module $\Z/\ell(-1)$, which is the  dual of
$\mu_\ell^{\otimes 2}$, has $G_F$-invariants iff $\Z/\ell\simeq
\mu_\ell$, i.e., iff $\mu_\ell\subset F$. Now we are done.

\bigskip

\section*{Appendix}

\medskip

In this Appendix, we will collect certain facts and results we use
concerning the flat topology, locally of finite type (cf.
\cite{Mi1}, chapter II).

\medskip

{\bf A1. \, The Setup}

\medskip

As in the main body of the paper, we assume that $F$ is a
non-archimedean {\em local field} of {\em characteristic zero}, with
{\em residue field} $k=F_q$ of {\em characteristic $p>0$}. Let
$S=\Spec\,\ScrO_F$, $j\colon\, U=\Spec\,F\to S$ and $i\colon\,
s=\Spec\,k\to S$ where $s$ is the closed point.

Let $E$ be an {\em elliptic curve} defined on $F$, and let
$A=E\times E$. Let $\ScrE$ be the N\'eron model of $E$ over $S$, so
we have in particular $E(F)=\ScrE(S)$ where $E(F)$ denotes the
$F$-rational points of $E$ and $\ScrE(S)$ are the sections of
$\ScrE$ over $S$. We have $E=j^*\ScrE$ and $i^*\ScrE$ as pullbacks,
and $i^*\ScrE$ is the reduction of $E$. For any $n \geq 1$, let
${\ScrE}[n]$ denote the kernel of multiplication by $n$ on $\ScrE$.

Let $\ell$ be a prime number $>2$, possibly equal to $p$. Under
these assumptions we can apply Proposition 2.5.1 and Corollary 2.5.2
to $E$ and $A=E\times E$. Since $H_{\et}^1(\Ebar,\Z_\ell(1))$ is the Tate module
${\mathcal T}_\ell(E)$, we have the homomorphism $c_\ell\colon\, T_F(A)\to
H^2(F,{\mathcal T}_\ell(E)^{\otimes 2})$. Reducing modulo $\ell$, we obtain now
a homomorphism
$$
s_\ell\colon\, T_F(A)/\ell\longrightarrow H^2(F,E[\ell]^{\otimes
2}).
$$

\medskip

{\bf A2. \, Flat topology and the integral part}

\medskip

We shall need flat topology, or to be precise, {\it flat topology,
locally of finite type} (cf. \cite{Mi1}).

Let $\ScrF$ be a sheaf on the {\it big flat site} $S_{\fl}$.
Consider the restriction $\ScrF_{\vert_U} = \ScrF_U:= j^\ast \ScrF$
to $U$. Now assume that $\ScrF_U$ is the pull back
$\pi_U^\ast(\ScrG)$ of a sheaf $\ScrG$ on the {\it \'etale site}
$U_{\et}$ under the morphism of sites $\pi_U: U_{\fl} \, \to \,
U_{\et}$. Then we have, in each degree $i$, a homomorphism
$$
\pi_U^\ast: H^i_{\et}(F, {\ScrG}\vert_U) \, \rightarrow \,
H^i_{\fl}(F, {\ScrF}\vert_U). \leqno(A.2.1)
$$
We will need the following well known fact (cf. \cite{Mi1}, Remark
3.11(b), pp. 116--117):

\medskip

\noindent{\bf Lemma A.2.2} \, \it This is an isomorphism when
$\ScrG$ is a locally constructible sheaf on $U_{\et}$. \rm

\medskip

Now let $\ScrF$ and $\ScrG$ be as above, with $\ScrG$ locally
constructible. Composing the natural map $H^i_{\fl}(S, \ScrF)
\rightarrow H^i_{\fl}(U, \ScrF)$ with the inverse of the isomorphism
of Lemma A.2.2, we get a homomorphism
$$
\beta: \, H^i(_\fl(S, \ScrF) \, \rightarrow \, H^i(F, \ScrG_F),
$$
where the group on the right is Galois cohomology in degree $i$ of
$F$, i.e., of Gal$(\overline F/F)$, with coefficients in the Galois
module $\ScrG_F$ associated to the sheaf $\ScrG$ on $U_\et$.

\medskip

\bigskip\noindent
{\bf Definition A.2.3} \, \it In the above setup, define the
integral part of $H^i(F, \ScrG_F)$ to be
$$
\overline H^i(F, {\ScrG_F}) : = \, {\text Im}(\beta).
$$
\rm

\medskip

The sheaves of interest to us below will be even constructible, and
in many cases, $\ScrG$ will itself be the restriction to $U$ of a
sheaf on $S_\et$. A key case will be when there is a finite, flat
groupscheme $C$ on $S$ and $\ScrG=C_U$, where $C_U$ is the sheaf on
$U_\et$ defined by $C_{\vert_U}$. We will, by abuse of notation, use
the same letter $\tilde C$ to indicate the corresponding sheaf on
$S_\fl$ or $S_\et$ as long as the context is clear.

\medskip

More generally, let $C_1, C_2, \dots, C_r$ denote a collection of
smooth, commutative, finite flat groupschemes over $S$, which are
\'etale over $U$. (They need not be \'etale over $S$ itself, and
possibly, $C_i=C_j$ for $i\ne j$.) Let $\ScrF_\et$, resp.
$\ScrF_{\fl}$, denote the sheaf for $S_{\et}$, resp. $S_{fl}$,
defined by the tensor product $\otimes_{j=1}^r \, \tilde C_j$. Since
each $\ScrG_i$ is \'etale over $U$, we have, for the restrictions to
$U$ via $\pi_U: U_{\fl} \to U_{\et}$, an isomorphism
$$
\ScrF_{\fl}\vert_U \, \simeq \, \ScrF_{\et}\vert_U. \leqno(A.2.4)
$$
This puts us in the situation above with $\ScrF=\ScrF_\fl$ and
$\ScrG=\ScrF_{\fl}\vert_U$.

\medskip

In fact, for the $\tilde C_i$ themselves, (A.2.4) follows from
\cite{Mi1}, p. 69, while for their tensor products, one proceeds via the
tensor product of sections of presheaves. Note also that
$j^\ast(\ScrF_\et)$ is constructible in $U_\et$ since it is locally
constant. Ditto for the kernels and cokernels of homomorphisms
between such sheaves on $U_\et$.

\medskip

{\bf A3. \, Application: \, The Image of Symbols}

\medskip

We begin with a basic result:

\medskip

\noindent {\bf Lemma A.3.1} \, \it We have for any prime $\ell$
(possibly $p$),
$$
\ScrO_F^*/\ell\arrowsim H_{\fl}^1(S,\mu_\ell)\arrowsim \overline H^1(F,
\mu_\ell)\hookrightarrow H_{\et}^1(F,\mu_\ell) \leftarrowsim F^*/\ell
$$
(with the natural inclusion). \rm

\medskip

\begin{Proof}
The first isomorphism follows from the exact sequence on ${S}_{\fl}$
$$
1\longrightarrow \mu_\ell\longrightarrow \ScrG_m{\overset
\ell\longrightarrow} \ScrG_m\longrightarrow 1
$$
and the fact that $H_{\fl}^1(S,\ScrG_m)=\Pic(S)=0$. (Of course, when
$\ell \ne p$, this sequence is also exact over $S_{\et}$. The second
isomorphism then follows from the inclusion $\ScrO_F^\ast/\ell
\hookrightarrow F^\ast/\ell$. The last two statements follow from
the definition or are well known.
\end{Proof}

\medskip

\noindent {\bf Proposition A.3.2} {\it Let $F$ be a local field of
characteristic $0$ with finite residue field $k$ of characteristic
$p$, $E/F$ an elliptic curve with good
reduction, and $\ScrE$ the N\'eron model. Then for any odd prime $\ell$,
possibly equal to $p$, there is a
short exact sequence
$$
0 \to E(F)/\ell \to H_{\fl}^1(S, \ScrE[\ell]) \to H^1(k,
\pi_0(\ScrE_{s}))[\ell] \to 0. \leqno(i)
$$
In particular, since $E$ has good reduction,
we have isomorphisms
$$
E(F)/\ell \simeq\ScrE(S)/\ell \arrowsim
H_{\fl}^1(S,\ScrE[\ell])\arrowsim \overline H^1(F,E[\ell])
$$}

\medskip

\begin{Proof}
Over $\mathcal{S}_{\fl}$ we have the following {\em exact sequence
of sheaves} associated to the group schemes (see \cite{Mi2}, p. 400,
Corollary C9):
$$
0\longrightarrow \ScrE[\ell]\longrightarrow \ScrE{\overset
\ell\longrightarrow}\ScrE,
$$
where the arrow on the right is surjective when $E$ has good
reduction. In any case, taking flat cohomology, we get an exact
sequence
$$
0 \to \ScrE(S)/\ell = E(F)/\ell \to H_{\fl}^1(S, \ScrE[\ell]) \to
H^1_\fl(S, \ScrE)[\ell].
$$

There is an isomorphism
$$
H_{\fl}^1(S,\ScrE) \, \simeq H^1(k,\pi_0(\ScrE_{s})).
\leqno(A.3.3)
$$
Since $E$ has good reduction, this is proved as follows:
The natural map
$$
H_{\fl}^1(S,\ScrE) \,
\rightarrow \, H^1(s,\ScrE_{s}), \, \, s = {\rm spec} k,
$$
is an isomorphism when $\ScrE$ is smooth over $S$ (see \cite{Mi1}, p.
116, Remark 3.11), and $H^1(s,\ScrE_{s})$ vanishes in this
case by a theorem of Lang (\cite{Lang1}.

Clearly, we now have the exact sequence (i).

\end{Proof}

\medskip

\noindent{\bf Proposition A.3.4} \, \it Let $E$ be an elliptic curve
over $F$ with good reduction, and let $\ell$ be an odd prime,
possibly the residual characteristic $p$ of $F$. Then the restriction $s_\ell^{\rm symb}$ of $s_\ell$
to the symbol group $ST_{F, \ell}(A)$ factors as follows:
$$
s_\ell^{\rm symb}: ST_{F, \ell}(A) \, \rightarrow \, \Hbar^2(F, E[\ell]^{\otimes
2}) \, \hookrightarrow \, H^2(F, E[\ell]^{\otimes 2}).
$$

\medskip

\begin{proof} Thanks to Proposition A.3.2 and the definitions, this
is a consequence of the following commutative diagram:
\begin{equation*}
\begin{CD}
E(F)/\ell \times E(F)/\ell &\rightarrow
&H_{\fl}^1(S,\ScrE[\ell]) \times H_{\fl}^1(S,\ScrE[\ell])&\\
@VVV \Big\downarrow{}\\
ST_{F, \ell}(A) &\,\longrightarrow \, &\,
H_{\fl}^2(S,\ScrE[\ell]^{\otimes 2}) &\longrightarrow &\, \overline
H^2(F,E[\ell]^{\otimes 2}) &\,\hookrightarrow
 & \, H^2(F,E[\ell]^{\otimes 2})
\end{CD}
\end{equation*}
\end{proof}

\medskip

{\bf A4. \, A lemma in the case of good ordinary reduction}

\medskip

Let $E$ be an elliptic curve over a non-archimedean local field $F$ of
characteristic zero, whose residue field $k$ has characteristic $p
>0$.

\medskip

\noindent{\bf Lemma A.4.1} \, Assume that $E$ has good ordinary
reduction, with N\'eron model $\ScrE$ over $S={\rm Spec}(\ScrO_F)$.
Then the following hold:
\begin{enumerate}
\item[(a)]There exists an exact sequence of group schemes, finite
and flat over $S$, and a corresponding exact sequence of sheaves on
$S_\fl$, both compatible with finite field extensions $K/F$:
$$
0\longrightarrow C\longrightarrow \ScrE[p]\longrightarrow C'
\longrightarrow 0,\leqno(A.4.2)
$$
where $C=\ScrE[p]^{\rm loc}$ and $C'=\ScrE[p]^{\rm et}$.
\item[(b)]If the points of $E_s^{\et}[p]$ are all rational
over $k$ then the above exact sequence becomes
$$
0\longrightarrow \mu_{p,S}\longrightarrow \ScrE[p]_S\longrightarrow
(\Z/p)_S \longrightarrow 0.\leqno(A.4.3)
$$
\item[(c)]If $E[p]$ is in addition semisimple as a Gal$(\overline F/F)$-module,
the exact sequence (A.4.3) splits as groupschemes (and {\it \`a
fortiori} as sheaves on $S_\fl$), i.e.,
$$
\ScrE[p] \, \simeq \, \mu_{p,S} \oplus (\Z/p)_S.\leqno(A.4.4)
$$
\end{enumerate}

\medskip

\begin{proof} \, (a) \, This sequence of finite, flat groupschemes is well known (see \cite{Tate4},
Thm. 3.4). The exactness of the sequence in $S_\fl$, which is clear
except on the right, follows from the fact that $\ScrE[p]\to
\ScrE[p]^\et$ is itself flat (cf. \cite{Raynaud1}, pp. 78--85).
Moreover, the naturality of the construction furnishes compatibility
with finite base change.

(b) \, Since $C'$ is $\ScrE[p]^\et$, it becomes $(\Z/p)_S$ when all
the $p$-torsion points of $\ScrE_s^\et$ become rational over $k$. On
the other hand, since $\ScrE[p]$ is selfdual, $C$ is the Cartier
dual of $C'$, and so when the latter is $(\Z/p)_S$, the former has
to be $\mu_{p, S}$.

(c) \, When $E[p]$ is semisimple, it splits as a direct sum $C_F
\oplus C'_F$. In other words, the groupscheme $\ScrE[p]$ splits over
the generic point, and since we are in the N\'eron model, any section
over the generic point furnishes a section over $S$. This leads to a
decomposition $C \oplus C'$
over $S$ as well. When we are in the situation of (b), the
semisimplicity assumption yields (A.4.4).
\end{proof}

\noindent{
\begin{tabular}{ll}
Jacob Murre &Dinakar Ramakrishnan\\
Department of Mathematics & Department of Mathematics\\
P.O.Box 9512, 2300 RA Leiden & 253-37 Caltech\\
University of Leiden & California Institute of Technology\\
Leiden, Netherlands & Pasadena, CA 91125, USA\\
murre(at)math.leidenuniv.nl & dinakar(at)caltech.edu
\end{tabular}}

\end{document}